\documentclass[12pt]{article}

\usepackage{amstext}
\usepackage{amsfonts}
\usepackage{amssymb}
\usepackage{amsbsy}
\usepackage{amsmath}

\newcounter{num} %

\newenvironment{theo}
{\refstepcounter{num}%
\bigskip\noindent{\bf Theorem~\arabic{num}. }\it}

\newenvironment{cor}
{\refstepcounter{num}%
\bigskip\noindent{\bf Corollary~\arabic{num}. }\it}

\newenvironment{lemma}
{\refstepcounter{num}%
\bigskip\noindent{\bf Lemma~\arabic{num}. }\it}

\newenvironment{prop}
{\refstepcounter{num}%
\bigskip\noindent{\bf Proposition~\arabic{num}. }\it}

\newtheorem{remark}{Remark}

\newcommand{\Ref}[1]{(\ref{#1})}

\newenvironment{proof}{\medskip\noindent{\it Proof. }}
{$\Box$ \bigskip}

\newenvironment{eq}{\begin{equation}}{\end{equation}}

\newcommand{\Pf}{{\rm pf}}
\renewcommand{\P}{{\rm{P }}}
\newcommand{\DP}{{\rm DP}}

\newcommand{\NN}{{\mathbb{N}} }
\newcommand{\ZZ}{{\mathbb{Z}} }
\newcommand{\QQ}{{\mathbb{Q}} }
\newcommand{\Sp}{S\!p}

\newcommand{\my}{}            
\newcommand{\myun}[1]{#1}
\newcommand{\My}{}            

\newcommand{\mmyun}[1]{\underline{#1}}
\newcommand{\MMy}[1]{#1}      

\newcommand{\Ring}{\mathcal{F}}                      
\newcommand{\Symmgr}{\mathcal{S}}                 
\newcommand{\Reprsp}{\mathcal{R}}                    
\newcommand{\Field}{\mathcal{K}}                    
\newcommand{\Quiver}{\mathcal{Q}}                    
\newcommand{\VectspV}{\mathcal{V}}                   
\newcommand{\VectspW}{\mathcal{W}}                   

\newcommand{\GroupG}{\mathcal{G}}                   
\newcommand{\GroupH}{\mathcal{H}}                   
\newcommand{\GroupL}{\mathcal{L}}                   
\newcommand{\GroupA}{\mathcal{A}}                   
\newcommand{\GroupB}{\mathcal{B}}                   
\newcommand{\GroupD}{\mathcal{D}}                   

\newcommand{\Dimvec}{\mmyun{d}}                   

\newcommand{\si}{\sigma}

\newcommand{\al}{\alpha}
\newcommand{\be}{\beta}

\newcommand{\ga}{\gamma}
\newcommand{\Ga}{\Gamma}
\newcommand{\de}{\delta}
\newcommand{\De}{\Delta}
\newcommand{\la}{\lambda}
\newcommand{\La}{\Lambda}

\newcommand{\LA}{\langle}
\newcommand{\RA}{\rangle}
\newcommand{\ov}[1]{\overline{#1}}

\newcommand{\tr}{{\rm tr}}

\newcommand{\diag}{\mathop{\rm diag}}

\newcommand{\sign}{\mathop{\rm{sgn }}}

\newcommand{\Hom}{{\mathop{\rm{Hom }}}}

\newcommand{\ot}{\mathop{\otimes}}

\begin{document}

 \title{Semi-invariants of mixed representations of quivers}
 \author{}

 \date{}



\maketitle 

\vspace{-1.5cm} 

$$
\begin{array}{cc}
\mbox{ \large A. A. Lopatin}&\mbox{\large A. N. Zubkov} \\
\mbox{ Institute of Mathematics ,}&\mbox{ Omsk State Pedagogical University,}\\
\mbox{ Siberian Branch of }&\mbox{Chair of Geometry,}\\
\mbox{ the Russian Academy of Sciences, }&\mbox{}\\
\mbox{ Pevtsova street, 13,}&\mbox{Tuhachevskogo embarkment, 14,}\\
\mbox{ Omsk 644099 Russia}&\mbox{Omsk 644099, Russia} \\
\mbox{ artem\underline{ }lopatin@yahoo.com}&\mbox{a.zubkov@yahoo.com} \\
\mbox{http://www.iitam.omsk.net.ru/\~{}lopatin/}&\\
\end{array}
$$

\bigskip

\begin{abstract}
A notion of a mixed representation of quivers can be derived from
ordinary quiver representation by considering the dual action of
groups on ``vertex"{} vector spaces together with their usual
action. A generating system for the algebra of semi-invariants of
mixed representations of a quiver is determined. This is done by
reducing the problem to the case of bipartite quivers of a special
form and using a function $\DP$ on three matrices, which is a
mixture of the determinant and two pfaffians.
\end{abstract}

\section{Introduction}\label{section_intro}

Quiver representations appeared for the first time
in~\cite{Gabriel72}. The importance of this notion is due to the
fact that the category of representations of a quiver is
equivalent to the category of finite-dimensional modules over the
path algebra associated with that quiver. Since any finite-dimensional basic algebra over an algebraically closed field is a
factor-algebra of the path algebra of some quiver (see Chapter~3
of~\cite{Kirichenko}), the set of finite-dimensional modules over
such an algebra is a full subcategory of the category of
representations of the quiver. A representation of a quiver
consists of a collection of vector spaces, assigned to its
vertices, and linear mappings between the vector spaces ``along"{}
the arrows. A morphism of two representations of the same quiver
is a collection of linear mappings between the vector spaces,
assigned to the same vertices, which commute with the linear
mappings of the representations. The set of quiver representations
of a fixed dimension can be endowed naturally with the structure
of a vector space and is called the space of quiver
representations of that dimension. Its group of automorphisms is a
direct product of the general linear groups, acting on ``vertex"{}
spaces. Obviously, the orbits of this action correspond to the
classes of isomorphic representations. If we compute generators of
the algebra of invariants, then we can distinguish closed orbits,
which correspond to semisimple representations.

At first this problem was solved in an important special case of a
quiver with one vertex and several loops. Its representation space
of dimension $n$ coincides with a space of several $n\times n$
matrices on which $GL(n)$ acts by diagonal conjugations. In the
case of a field of characteristic zero, generators for the algebra
can easily be found by means of the classical invariant theory
(for details see, e.g.,~\cite{Vinberg_Popov}, page~257). Defining
relations between generators were described in~\cite{Razmyslov74}
and~\cite{Procesi76}. In the case of positive characteristic,
generators for the algebra of invariants were described
in~\cite{Donkin92a}, and relations between them were established
in~\cite{Zubkov96}. Both of these results are formulated in terms
of coefficients $\si_k$ in the characteristic polynomial of a
matrix. Their proofs rely on the theory of modules with good
filtration~\cite{Donkin85}.
Subsequently methods from~\cite{Donkin92a} and~\cite{Zubkov96}
were successfully applied to representations of any
quiver~\cite{Donkin94} and~\cite{Zubkov_Fund_Math_01}. In the case
of a field of characteristic zero these results were proved
in~\cite{Le_Bruyn_Procesi_90} and later in~\cite{Domokos98}.

Natural generalizations of the above construction are invariants of quivers under the
action of classical groups, i.e., of $GL(\VectspV)$, $O(\VectspV)$, $\Sp(\VectspV)$,
$SL(\VectspV)$, and $SO(\VectspV)$. In other words, we assume that some classical group
acts on the vector space assigned to a vertex, so the product of these groups acts on the
space of quiver representations. In particular, the algebra of invariants under the
action of a product of the special linear groups is called the algebra of {\it
semi-invariants}. Its generators were calculated in~\cite{DZ01} using methods
from~\cite{Donkin92a},~\cite{Donkin94},~\cite{Zubkov96}, and~\cite{Zubkov_Fund_Math_01},
and, independently, in~\cite{DW00} and~\cite{DW_LR_02} utilizing methods of the
representation theory of quivers. In the paper~\cite{Schofield_van_den_Bergh_01}
generators for semi-invariants were given in the case of a field of characteristic zero.

For remaining classical groups, the first step was performed in
the classical work of Procesi~\cite{Procesi76}. He calculated
generators and relations between them for the orthogonal and
symplectic invariants of several matrices over a field of
characteristic zero. The study of the case of zero characteristic
was continued in~\cite{Aslaksen95}, where generators for the case
of the special orthogonal group were found. The results concerning
generators for the orthogonal and symplectic invariants were
generalized in~\cite{Zubkov99} to a field of positive
characteristic (assumed to be odd in the case of the orthogonal
group). The main reduction of paper~\cite{Zubkov99} (also
see~\cite{ZubkovI}) actually shows that the most general concept
for dealing with product of arbitrary classical groups is the
notion of mixed representations of quivers, introduced
in~\cite{ZubkovI}.

The main idea of the definition of {\it mixed representations} is
to deal with bilinear forms in addition to linear mappings. It
works as follows. Bilinear forms on some vector space $\VectspV$
are in one to one correspondence with linear mappings from
$\VectspV$ to its dual vector space $\VectspV^{\ast}$. The
standard action of $GL(\VectspV)$ on $\VectspV$ induces the action
on $\VectspV^{\ast}$, and, therefore, on
$\Hom(\VectspV,\VectspV^{\ast})$. So, together with the action of
the general linear groups on the vector spaces assigned to
vertices, we consider their dual action.

As an example, consider a
quiver $\Quiver$ that can schematically be depicted as %
$$\bullet{v} \stackrel{a,b}{\Longrightarrow} \bullet{w}\quad, $$
where $v,w$ are vertices; $a,b$ are arrows from $v$ to $w$; a
vector space $\VectspV$ is assigned to $v$, and $\VectspV^{\ast}$
is assigned to $w$. Then the orbits of $GL(\VectspV)$ on the space
of mixed representations of $\Quiver$ correspond to pairs of
bilinear forms on $\VectspV$. The classification problem for pairs
of symmetric and skew-symmetric bilinear forms is a classical
topic going back to Weierstrass and Kronecker
(see~\cite{Gantmacher},~\cite{Hodge_Pedoe}, and~\cite{Gurevich50}).

Another motivation for introducing mixed representations are
papers~\cite{DW02} and~\cite{Shmelkin}, where {\it orthogonal} and
{\it symplectic} representations of {\it symmetric} quivers and
{\it symmetric} representations of {\it signed} quivers,
respectively, were introduced. These notions relate to mixed
representations of quivers very closely. In fact, such
representations are mixed representations for some products of
classical linear groups. In~\cite{DW02} it was established that
representations of a symmetric quiver are in one-to-one
correspondence with a some subset of representations of the
associated double quiver. Consequently the symmetric quivers of
tame and finite type were classified. These results were
generalized to signed quivers in~\cite{Shmelkin}.

In~\cite{ZubkovI} generating invariants for mixed representations
of quivers were found, as well as invariants for representations
of symmetric and signed quivers. In a subsequent
paper~\cite{ZubkovII} defining relations for the algebra of
invariants of mixed representations of quivers were established.
To describe relations, together with $\si_k$ a new function
$\si_{k,r}$ on three matrices is used. This function describes
some relations
for the quiver schematically depicted as %
$$c\subset \bullet{v} \stackrel{a,b}{\Longleftrightarrow} \bullet{w} $$
Here the arrow $a$ goes from $v$ to $w$, the arrow $b$ goes in the
opposite direction, $c$ is a loop on the vertex $v$, and the
vector space assigned to $w$ is dual to the space assigned to $v$.
The general results were applied to obtain information about
defining relations for the orthogonal and symplectic invariants
(see Sections~3 and~4 of~\cite{ZubkovII}).

The current paper is dedicated to semi-invariants of mixed
representations of quivers. In general we follow the approach
of~\cite{DZ01}. The main result of the paper is a description of a
generating system for semi-invariants of mixed representations of
quivers.

Section~\ref{section_2} contains necessary definitions and some
auxiliary results.

In Section~\ref{section_F} the function $\DP$ on three matrices is
introduced and some of its combinatorial properties are studied.
This function is a mixture of the determinant and two pfaffians.
Note that $\DP$ relates to $\si_{k,r}$ in the same way as the
determinant relates to $\si_k$, i.e., for $1\leq k\leq t$ the
function $\si_{k,r}(X,Y,Z)$ is the coefficient of $\la^{t-k}$ in
the polynomial $\DP_{r,r}(X+\la E,Y,Z)$ for $(t+2r)\times(t+2r)$
matrices $X,Y,Z$. (This is a direct consequence of the
decomposition formula from~\cite{Lop_bplp}.)

In Section~\ref{section_4} the general problem is reduced to the
case of the so-called zigzag-quivers, which are special cases of
bipartite quivers. We prove that semi-invariants of mixed
representations of a quiver are spanned by semi-invariants of
mixed representations of some zigzag-quiver and describe
explicitly this reduction (see Theorem~\ref{theo_reduction}).

In Sections~\ref{section_results} and~\ref{section_proof} we
determine a generating system for semi-invariants of mixed
representations of a zigzag-quiver. There is an obvious way to
obtain some semi-invariants. First we consider triplets of block
matrices, using generic matrices of mixed representations of a
quiver as blocks. Partial linearizations of the function $\DP$ on
these triplets are semi-invariants of mixed representations of a
quiver. We prove that all semi-invariants of mixed representations
of a quiver belong to the linear span of these semi-invariants
(see Theorem~\ref{theo1}). In Subsection~\ref{results_special} a
special case is considered. Theorem~\ref{theo1} is proved in
Section~\ref{section_proof}.

The paper is concluded by example given in
Section~\ref{section_appl}.

Using the main result of this paper, together with a reduction to
mixed representations of quivers (cf.~\cite{ZubkovI}) and the
decomposition formula from~\cite{Lop_bplp}, the first author
completed a description of generators for the invariants under the
action of a product of classical groups on the space of (mixed)
representations of a quiver in~\cite{Lop_so_inv}. In particular,
generators for the invariants of several matrices under the action
of the special orthogonal group were found, where the
characteristic is different from $2$ when we consider the
(special) orthogonal group.

\section{Preliminaries}\label{section_2}

\subsection{Notations and remarks}\label{section_2.1}

Let $\Field$ be an infinite field of arbitrary characteristic. All
vector spaces, algebras and modules are over $\Field$ unless
otherwise stated. Denote by ${\NN}$ the set of all non-negative
integers, by $\ZZ$ the set of all integers, and by $\QQ$ the
quotient field of $\ZZ$. We use capital letters like
$\mathcal{A}$, $\mathcal{B}$, $\mathcal{C}$,  etc.~for sets
endowed with some algebraic structure and for quivers (see
Section~\ref{section_repr_quiv}).

For a vector space $\VectspV$ let $\VectspV^{\ast}$,
$S(\VectspV)$, $S^t(\VectspV)$, $\otimes^t \VectspV$, and
$\wedge^t \VectspV$, respectively, stand for the dual space, the
symmetric algebra, the $t$-th symmetric power, the $t$-th tensor
power, and the $t$-th exterior power of $\VectspV$, respectively.

Suppose a reductive algebraic group $\GroupG$ acts on vector
spaces $\VectspV$ and $\VectspW$, and let $v\in \VectspV$, $w\in
\VectspW$, $g\in \GroupG$. This action induces an action on
\begin{itemize}
\item the coordinate algebra $\Field[\VectspV]$ of the affine
variety $\VectspV$: for $f\in \Field[\VectspV]$ we put $(g\cdot
f)(v)=f(g^{-1}\cdot v)$;

\item the dual space $\VectspV^{\ast}$, which we consider as the
degree one homogeneous component of the graded algebra
$\Field[\VectspV]$;

\item  the space $\Hom_{\Field}(\VectspV,\VectspW)$ of
$\Field$-linear maps from $\VectspV$ to $\VectspW$: for $H\in
\Hom_{\Field}(\VectspV,\VectspW)$ put $(g\cdot H)(v)=g\cdot
(H(g^{-1}\cdot v))$;

\item the tensor product $\VectspV\otimes \VectspW$: $g\cdot
(v\otimes w)=g\cdot v\otimes g\cdot w$.
\end{itemize}

Given $n\in\NN$, consider the vector space $\VectspV(n)=\Field^n$
of column vectors of length $n$ with the {\it standard} basis
$v_1,\ldots,v_n$, where $v_i$ is a column vector, whose $i$-th
entry is $1$ and the rest of entries are zeros. Consider the dual
space $\VectspV(n)^{\ast}$ with the dual basis $v_1^{\ast},\ldots,
v_n^{\ast}$. We identify $\VectspV(n)^{\ast}$ with the space of
column vectors of length $n$, so $v_i^{\ast}$ is the same column
vector as $v_i$. Denote by $\Field^{n\times m}$ the space of
$n\times m$ matrices.

The general linear group $GL(n)=\{g\in \Field^{n\times
n}\,|\,\det(g)\neq0\}$ acts on $\VectspV(n)$ by left
multiplication. Let $\{1, \ast\}$ be the group of order two, where
the symbol $1$ stands for the identity element of the group. Given
$\al\in\{1, \ast\}$, $g\in GL(n)$, define
$$\VectspV(n)^{\al}=\left\{
\begin{array}{cc}
\VectspV(n)&\!\!\!\!\text{for } \al=1\\
\VectspV(n)^{\ast}&\!\!\!\!\text{for } \al=\ast\\
\end{array}
\right.,\quad v_i^{\al}=\left\{
\begin{array}{cc}
v_i&\!\!\!\!\text{for } \al=1\\
v_i^{\ast}&\!\!\!\!\text{for } \al=\ast\\
\end{array}
\right.,\quad g^{\al}=\left\{
\begin{array}{cc}
g&\!\!\!\!\text{for } \al=1\\
(g^{-1})^{T}&\!\!\!\!\text{for } \al=\ast\\
\end{array}
\right..$$ %
Note that
\begin{eq}\label{eq_action}
g\cdot v_i^{\al}=g^{\al}v_i^{\al}.
\end{eq}

Let $\al,\be\in\{1,\ast\}$, $\VectspV=\VectspV(n)^{\al}$,
$\VectspW=\VectspW(m)^{\be}$ and let $w_1,\ldots,w_m$ be the
standard basis for $\VectspW(m)$. Let groups $\GroupG_1\subset
GL(n)$, $\GroupG_2\subset GL(m)$ be closed under transposition and
the group $\GroupG=\GroupG_1\times \GroupG_2$ acts on $\VectspV$,
$\VectspW$ as follows: for $g=(g_1,g_2)\in \GroupG$, $v\in
\VectspV$, $w\in \VectspW$ we have $g\cdot v=g_1\cdot
v=g_1^{\al}v$, $g\cdot w=g_2\cdot w=g_2^{\be}w$. Let
$\Field[\Hom_{\Field}(\VectspV,\VectspW)]=\Field[x_{ij}\,|\,1\leq
i\leq m,\,1\leq j\leq n]$ be the coordinate algebra of the space
$\Hom_{\Field}(\VectspV,\VectspW)$, and $X=(x_{ij})$ be the
$m\times n$ matrix. Here $x_{ij}$ maps a matrix $H\in
\Field^{m\times n}\simeq \Hom_{\Field}(\VectspV,\VectspW)$ to the
$(i,j)$-th entry of $H$, where the isomorphism is determined by
the choice of bases for $\VectspV,\VectspW$. Consider
$g=(g_1,g_2)\in \GroupG$, $H\in \Field^{m\times n}
\simeq\Hom_{\Field}(\VectspV,\VectspW)$. The following lemma is
proved by straightforward calculations.

\begin{lemma}\label{lemma1}
Using the preceding notation we have:

1) $g\cdot H=g_2^{\be}H(g_1^{\al})^{-1}$.

2) $g\cdot X=(g_2^{\be})^{-1}Xg_1^{\al}$, where $g\cdot X$ means
the matrix whose $(i,j)$-th entry is equal to $g\cdot x_{ij}$.

3) The homomorphism of spaces
$\Phi:\Field[\Hom_{\Field}(\VectspV,\VectspW)]\to
S(\VectspW(m)^{\be\ast} \otimes \VectspV(n)^{\al})$, defined by
$\Phi(x_{ij})=w_i^{\be\ast} \otimes v_j^{\al}$, is an isomorphism
of $\GroupG$-modules.
\end{lemma}
\begin{proof}
1) Let $v\in \VectspV$. Then $(g\cdot H)(v)=g\cdot (H(g^{-1}\cdot
v))=g\cdot(H(g_1^{-1})^{\al}v)= g_2^{\be}H (g_{1}^{-1})^{\al} v$.

2) For $f\in \Field[x_{ij}]$ and $H\in \Field^{m\times n}\simeq
\Hom_{\Field}(\VectspV,\VectspW)$ we have $(g\cdot
f)(H)=f(g^{-1}\cdot H)=f((g_{2}^{\be})^{-1}Hg_1^{\al})$.

3) We prove that $\Phi$ is $\GroupG$-homomorphism. Denote the
entries of the matrix $g$ by $g_{ij}$. By the previous item,
$\Phi(g\cdot x_{ij})$ equals
$$\Phi(((g_2^{\be})^{-1}X g_1^{\al})_{ij})=
\Phi(\sum\nolimits_{r,s}((g_2^{\be})^{-1})_{ir}x_{rs}(g_1^{\al})_{sj})=
\sum\nolimits_{r,s}((g_2^{\be})^{-1})_{ir} (g_1^{\al})_{sj}
w_r^{\be\ast}\otimes v_s^{\al}.$$ On the other hand, $g\cdot
(w_i^{\be\ast}\otimes v_j^{\al})= g\cdot w_i^{\be\ast}\otimes
g\cdot v_j^{\al}= g_2^{\be\ast}w_i^{\be\ast}\otimes
g_1^{\al}v_j^{\al}= \sum_{r,s}(g_2^{\be\ast})_{ri}
(g_1^{\al})_{sj} w_r^{\be\ast}\otimes v_s^{\al}.$ The required
statement follows from $(g_2^{\be\ast})^T=(g_2^{\be})^{-1}$.
\end{proof}

Denote by $id_n$ the identity permutation of $\Symmgr_n$. Given
permutations $\si_1\in \Symmgr_{n_1},\ldots,\si_k\in
\Symmgr_{n_k}$, denote the product of their signs
$\sign(\si_1)\cdot\ldots\cdot\sign(\si_k)$ by
$\sign(\si_1\ldots\si_k)$. For positive integers $n_1,n_2,n_3$ and
permutations $\si_1\in \Symmgr_{n_1},\si_2\in
\Symmgr_{n_2},\si_3\in \Symmgr_{n_3}$ we regard $\si=\si_1\times
\si_2\times \si_3$ as an element of $\Symmgr_{n_1+n_2+n_3}$ given
by
$$
\si(l)=\left\{
\begin{array}{ccl}
\si_1(l),&{\rm if}& l\leq n_1\\
\si_2(l-n_1)+n_1,&{\rm if}& n_1< l\leq n_1+n_2\\
\si_3(l-n_1-n_2)+n_1+n_2,&{\rm if}& n_1+n_2< l\leq n_1+n_2+n_3.\\
\end{array}
\right.
$$

For any vector spaces $\VectspV,\VectspW$ and a positive integer
$n$ there is a homomorphism of vector spaces
$$\zeta_n:\wedge^n \VectspV^{\ast}\otimes\wedge^n \VectspW\to S^{n}(\VectspV^{\ast}\otimes
\VectspW),$$ %
defined by
$$\zeta_n(v_1^{\ast}\wedge\ldots\wedge
v_n^{\ast} \otimes w_1\wedge\ldots\wedge
w_n)=\sum\nolimits_{\rho\in
\Symmgr_n}\sign(\rho)\prod\nolimits_{l=1}^n v_l^{\ast}\otimes
w_{\rho(l)},
$$ %
where $v_i^{\ast}\in \VectspV^{\ast},w_j\in \VectspW$. Consider an
ordered set $\myun{\ga}=(\my\ga_1,\ldots,\my\ga_p)\in\NN^p$. We
set
$\wedge^{\myun{\ga}}\VectspV=\wedge^{\my\ga_1}\VectspV\otimes\ldots\otimes\wedge^{\my\ga_p}\VectspV$.
Define a homomorphism of vector spaces
$$\zeta_{\myun{\ga}}:
\wedge^{\myun{\ga}}\VectspV^{\ast}\otimes\wedge^{\myun{\ga}}\VectspW\to
S^{\my\ga_1+\ldots+\my\ga_p}(\VectspV^{\ast}\otimes \VectspW)$$
by $\zeta_{\myun{\ga}}(c)=\zeta_{\my\ga_1}(a_1\otimes b_1)\ldots
\zeta_{\my\ga_p}(a_p\otimes b_p)$, where
$c=(a_1\otimes\ldots\otimes a_p)\otimes (b_1\otimes\ldots\otimes
b_p)$, and  $a_i\in \wedge^{\my\ga_i}\VectspV^{\ast}$, $b_i\in
\wedge^{\my\ga_i}\VectspW$ ($1\leq i\leq p$).

The symbol $a\to b$ will indicate a substitution. For $a\in\NN$
and $B\subseteq\NN$ we write $a+B$ for $\{a+b\,|\,b\in B\}$.
Denote the Kronecker symbol by $\de_{ij}$. Given a group $\GroupG$
and $a,b\in \GroupG$, we set $a^b=b^{-1}ab$. For integers $i<j$
denote the interval $i,i+1,\ldots,j-1,j$ by $[i,j]$. For the
cardinality of a set $A$ we write $\#A$.

\subsection{The space of mixed representations of a
quiver}\label{section_repr_quiv}

A {\it quiver} $\Quiver=(\Quiver_0,\Quiver_1)$ is a finite
oriented graph, where $\Quiver_0=\{1,\ldots,l\}$ is the set of
vertices and $\Quiver_1$ is the set of arrows. For an arrow $a$,
denote by $a''$ its tail, and by $a'$ its head. Consider finite-dimensional vector spaces $\VectspV_1=\Field^{n_1},\ldots,
\VectspV_l=\Field^{n_l}$ of column vectors and fix their standard
bases. The vector $\mmyun{n}=(n_1,\ldots,n_l)$ is called the {\it
dimension} vector. Consider $\myun{\al}=(\al_1,\ldots,\al_l)$,
where $\al_1,\ldots,\al_l\in\{1,\ast\}$ and assign the vector
space $\VectspV_u^{\al_u}$ to each $u\in \Quiver_0$.

Further, assume that there is an involution $\phi: \Quiver_0\to
\Quiver_0$ satisfying the properties:
\begin{enumerate}
\item[a)] $n_{u}=n_{\phi(u)}$;

\item[b)] $\phi(u)=u$ implies $\al_u=1$;

\item[c)] $\phi(u)\neq u$ implies $\al_{\phi(u)}=\ast\al_u$.
\end{enumerate}
In other words, some of the vertices are grouped into
nonintersecting pairs, where vector spaces corresponding to a
pair of vertices have equal dimensions and one of the vector
spaces from such a pair is dual to another whereas all the
remaining vertices $u$ have $\al_u = 1$ and are stable under
$\phi$.

Define $GL(\mmyun{n})=GL(n_1)\times\ldots\times GL(n_l)$ and  %
$$\Reprsp=\Reprsp_{\myun{\al},\phi}(\Quiver,\mmyun{n})=%
\bigoplus\limits_{a\in \Quiver_1}\Field^{n_{a'}\times n_{a''}}\simeq%
\bigoplus\limits_{a\in \Quiver_1}\Hom_{\Field}(\VectspV_{a''}^{\al_{a''}},\VectspV_{a'}^{\al_{a'}}),$$ %
where the isomorphism is given by the choice of bases for
$\VectspV_1,\ldots,\VectspV_l$. For any vertex $u$ vector spaces
$\VectspV_u$ and $\VectspV_u^{\ast}$ are $GL(n_u)$-modules. Then
the group $GL(\mmyun{n})$ acts on $\Reprsp$ by the rule: for
$g=(g_1,\ldots,g_l)\in GL(\mmyun{n})$, $(H_a)_{a\in \Quiver_1}\in
\Reprsp$ we have
$$g\cdot (H_a)_{a\in \Quiver_1}=(g\cdot H_a)_{a\in \Quiver_1}=
(g_{a'}^{\al_{a'}}H_a(g_{a''}^{\al_{a''}})^{-1})_{a\in
\Quiver_1}$$ (see Lemma~\ref{lemma1}). For each pair $(u,\phi(u))$
with $u\in \Quiver_0$ and $\phi(u)\neq u$ replace the factor
$GL(n_{u})\times GL(n_{\phi(u)})$ of the group $GL(\mmyun{n})$ by
its diagonal subgroup. The resulting subgroup
$$\{g\in GL(\mmyun{n})\,|\,g_u=g_{\phi(u)}\text{ for all }u\in \Quiver_0\}$$
is isomorphic to
$$GL_{\myun{\al},\phi}(\mmyun{n})=\prod_{a\in \Quiver_0,\,\al(a)=1} GL(n_{a})$$
by $(g_a)_{a\in \Quiver_0}\mapsto (g_a)_{a\in
\Quiver_0,\,\al(a)=1}$. This isomorphism induces the action of
$GL_{\myun{\al},\phi}(\mmyun{n})$ on $\Reprsp$. So, one and the
same general linear group acts  on both $\VectspV_u$ and
$\VectspV_u^{\ast}$. The space $\Reprsp$ together with the action
of $GL_{\myun{\al},\phi}(\mmyun{n})$ on it is called an 
$\mmyun{n}$-dimensional {\it space of mixed representations} of
$\Quiver$, and elements of $\Reprsp$ are called {\it mixed
representations}.

The coordinate ring of the affine variety $\Reprsp$ is the polynomial ring  %
$$\Field[\Reprsp]=\Field[x_{ij}^b\,|\,b\in \Quiver_1,\,1\leq i\leq n_{b'},1\leq j\leq n_{b''}].$$ %
Here $x_{ij}^b$ stands for the coordinate function on $\Reprsp$
that takes a representation $H=(H_a)_{a\in \Quiver_1}$ to the
$(i,j)$-th entry of  $H_{b}$. The ring  $\Field[\Reprsp]$ is a
$GL_{\myun{\al},\phi}(\mmyun{n})$-module
(Section~\ref{section_2.1}).
Denote by %
$$\Field[\Reprsp]^{GL_{\myun{\al},\phi}(\mmyun{n})}=\{f\in \Field[\Reprsp]\,|\,g\cdot f=f \mbox{ for all }
g\in GL_{\myun{\al},\phi}(\mmyun{n})\} $$ %
the algebra of {\it invariants} of mixed representations of the
quiver $\Quiver$.

Let
$SL_{\myun{\al},\phi}(\mmyun{n})=GL_{\myun{\al},\phi}(\mmyun{n})\cap
(SL(n_1)\times \ldots \times SL(n_l))$.
Denote by %
$$\Field[\Reprsp]^{SL_{\myun{\al},\phi}(\mmyun{n})}=\{f\in \Field[\Reprsp]\,|\,g\cdot f=f \mbox{ for all }
g\in SL_{\myun{\al},\phi}(\mmyun{n})\} $$ %
the algebra of {\it semi-invariants} of mixed representations of
the quiver.

If $\al_u=1$ for all $u\in \Quiver_0$, then (semi)-invariants of
mixed representations of $\Quiver$ are (semi)-invariants of
representations of $\Quiver$.

For $\myun{\epsilon}=(\epsilon_1,\ldots,\epsilon_m)\in\ZZ^m$,
where $m$ is the cardinality of $\{u\in \Quiver_0\,|\,\al_u=1\}$,
denote the space of {\it relative invariants of  weight}
$\myun{\epsilon}$ by
$$\Field[\Reprsp]^{GL_{\myun{\al},\phi}(\mmyun{n}),\myun{\epsilon}}=\{f\in \Field[\Reprsp]\,|\,g\cdot f=
(\prod_{u\in \Quiver_0,\,\al_u=1}{\det}^{\epsilon_u}(g_u))f \mbox{
for all }
g\in GL_{\myun{\al},\phi}(\mmyun{n})\}. $$ %

\subsection{Distributions, partitions and Young subgroups}\label{section_Young}

By a {\it distribution} $B=(B_1,\ldots,B_d)$ of a set $[1,t]$ we
mean an ordered partition of the set $[1,t]$ into pairwise
disjoint subsets $B_j$ ($1\leq j\leq d$), which are called
components of the distribution. To every $B$ we associate two
functions $l\mapsto B|l|$ and $l\mapsto B\LA l\RA$ ($1\leq l\leq
t$), defined by the rules:
$$B|l|=i, \text{ if }l\in B_i,\text{ and }B\LA
l\RA=\#\{[1,l]\cap B_i\,|\,l\in B_i\}.$$ %
The symmetric group $\Symmgr_t$ acts on $[1,t]$ and contains the
{\it Young
subgroup} %
$$\Symmgr_B=\{\pi\in \Symmgr_t\,|\, \pi(B_i)=B_i \text{ for }1\leq i\leq
d\}.$$ %
For $\si\in \Symmgr_t$ we set
$B^{\si}=(B_1^{\si},\ldots,B_d^{\si})$, where
$B_i^{\si}=\si^{-1}B_i$, $1\leq i\leq d$. The intersection $A\cap
B$ of distributions $A=(A_1,\ldots,A_p)$ and $B=(B_1,\ldots,B_q)$
of the same set is the result of pairwise intersections of all
components of the given distributions, i.e., $C=(C_1,\ldots,C_d)$,
where for any $1\leq k\leq d$ we have $C_k\neq\emptyset$ and there
exist $i,j$ such that $1\leq i\leq p$, $1\leq j\leq q$, and
$C_k=A_i\cap B_j$; moreover, we assume that for any $1\leq
k_1<k_2\leq d$ the minimum of $C_{k_1}$ is less that the minimum
of $C_{k_2}$. If $A,B$ are distributions of the same set, then we
write $A\leq B$ provided each component of $A$ is contained in
some component of $B$.

A vector $\mmyun{t}=(t_1,\ldots,t_d)\in \NN^d$ determines the
distribution $\MMy{T}=(\MMy{T}_1,\ldots,\MMy{T}_d)$ of the set
$[1,t]$, where $t=t_1+\ldots +t_d$ and
$\MMy{T}_i=\{t_1+\ldots+t_{i-1}+1,\ldots,t_1+\ldots+t_i\}$, $1\leq
i\leq d$. We use capital letters to refer to the distribution
determined by a vector. Note that for $1\leq l\leq t$ we have
$$\MMy{T}|l|=i,\text{ if }t_1+\ldots+t_{i-1}<l\leq t_1+\ldots+t_{i}, \text{ and}$$
$$\MMy{T}\LA l\RA=l-(t_1+\ldots+t_{i-1}),\text{ where }\MMy{T}|l|=i.$$

It is easy to derive the following lemma (see Section~2.2
of~\cite{DZ01}).

\begin{lemma}\label{lemma2}
Given a distribution $B=(B_1,\ldots,B_d)$ of a set $[1,t]$,
consider $\mmyun{t}=(\#B_1,\ldots,\#B_d)$. Let $\MMy{T}$ be the
distribution determined by $\mmyun{t}$. Then for each $\si\in
\Symmgr_t$ we have:

1. $\Symmgr_{B^{\si}}=\si^{-1}\Symmgr_B\si=\Symmgr_B^{\si}$,
$B^{\si}|l|=B|\si(l)|$ for any $1\leq l\leq t$.

2. There is a permutation $\eta\in \Symmgr_B^{\si}$ with $B\LA
\si(l)\RA=B^{\si}\LA \eta(l)\RA$ for any $1\leq l\leq t$.

3. There exists a permutation $\rho\in \Symmgr_t$ with
$B=\MMy{T}^{\si}$ for any $\si\in \Symmgr_{\MMy{T}}\rho$.
Moreover, there is a unique permutation $\si\in
\Symmgr_{\MMy{T}}\rho$ that satisfies  the condition
$\MMy{T}^{\si}\LA l\RA=\MMy{T}\LA\si(l)\RA$ for any $1\leq l\leq
t$.
\end{lemma}
\bigskip

A vector $\myun{\ga}=(\my\ga_1,\ldots,\my\ga_p)\in\NN^p$
satisfying $\my\ga_1\geq\ldots \geq\my\ga_p$ and
$\my\ga_1+\ldots+\my\ga_p=t$ is called a {\it partition} of
$t\in\NN$ and is denoted by $\myun{\ga}\vdash t$. A {\it
multipartition} $\myun{\ga}\vdash \mmyun{t}$ is a $d$-tuple of
partitions $\myun{\ga}=(\myun{\ga}(1),\ldots,\myun{\ga}(d))$,
where $\myun{\ga}(i)\vdash t_i$, $\mmyun{t}=(t_1,\ldots,t_d)$. We
identify a multipartition $\myun{\ga}$ with a vector
$(\my\ga_1,\ldots,\my\ga_p)\in\NN^p$, where
$\myun{\ga}(i)=(\my\ga_{p_1+\ldots +p_{i-1}+1},\ldots,
\my\ga_{p_1+\ldots+p_i})$, $p=p_1+\ldots +p_d$, and say that
$\mmyun{p}=(p_1,\ldots,p_d)$ is the {\it height} of $\myun{\ga}$.
Any multipartition $\myun{\ga}\in\NN^{p}$ determines the
distribution $\My\Ga=(\My\Ga_1,\ldots,\My\Ga_p)$ and the Young
subgroup $\Symmgr_{\My\Ga}$, as above.

For a subgroup $\GroupA$ of $\GroupG$ denote by $\GroupG/\GroupA$,
$\GroupA\backslash \GroupG$ and $\GroupA\backslash \GroupG/B$,
respectively, some system of representatives for the left, right
and double cosets, respectively. Note that if $\si$ ranges over
$\GroupG/\GroupA$, then $\si^{-1}$ ranges over $\GroupA\backslash
\GroupG$ and vice versa.

\begin{lemma}\label{lemma4}
Let $\GroupG$ be a group, and let $\GroupL,\GroupA,\GroupB$ be its
subgroups. Then
\begin{enumerate}
\item[a) ] for any system of representatives $\GroupL/ \GroupL\cap
\GroupB$ there exist systems of representatives $\GroupL\cap
\GroupA\backslash \GroupL/\GroupL \cap \GroupB$ and $\GroupA\cap
\GroupL/\GroupA\cap \GroupL\cap \GroupB^{\pi}$, where $\pi^{-1}\in
\GroupL\cap
\GroupA\backslash \GroupL/\GroupL \cap \GroupB$, such that %
\begin{eq}\label{eq_systems_of_representatives}
\GroupL/ \GroupL\cap \GroupB=\{\nu\pi^{-1}\,|\,\pi^{-1}\in
\GroupL\cap \GroupA\backslash \GroupL/\GroupL \cap \GroupB ,\,
\nu\in \GroupA\cap \GroupL/\GroupA\cap \GroupL\cap
\GroupB^{\pi}\};
\end{eq} %
\item[b) ] for any systems of representatives $\GroupL\cap
\GroupA\backslash \GroupL/\GroupL \cap \GroupB$ and $\GroupA\cap
\GroupL/\GroupA\cap \GroupL\cap \GroupB^{\pi}$, where $\pi^{-1}\in
\GroupL\cap \GroupA\backslash \GroupL/\GroupL \cap \GroupB$, there
exist a system of representatives $\GroupL/ \GroupL\cap \GroupB$
such that~\Ref{eq_systems_of_representatives} holds.
\end{enumerate}
\end{lemma}
\begin{proof}
\textbf{ a)} Consider a set  $\GroupL/ \GroupL\cap
\GroupB=\{g_i\,|\,i\in I\}$. Since $\GroupL$ is a union of double
cosets,  for each $i\in I$ there is a $\pi_i\in \GroupL$ such that
$g_i\in (\GroupL\cap \GroupA)\pi_i^{-1}(\GroupL\cap \GroupB)$.
Multiplying $g_i$ by an appropriate element from $\GroupL\cap
\GroupB$ we can assume that $g_i=\nu\pi_i^{-1}$ for some $\nu\in
\GroupL\cap \GroupA$ for all $i$. If two elements $\nu\pi^{-1}$
and $\nu'\pi'^{-1}$ from $\{\nu\pi^{-1}\,|\, \pi^{-1}\in
\GroupL\cap \GroupA\backslash \GroupL/\GroupL \cap \GroupB, \nu\in
\GroupA\cap \GroupL\}$ are equal modulo $\GroupL\cap \GroupB$,
then $(\GroupL\cap \GroupA)\pi^{-1}(\GroupL\cap
\GroupB)=(\GroupL\cap \GroupA)\pi'^{-1}(\GroupL\cap \GroupB)$, and
thus $\pi=\pi'$. Therefore $\nu\pi^{-1}\equiv\nu'\pi'^{-1}\;{\rm
mod}(\GroupL\cap \GroupB)$ holds if and only if $\nu^{-1}\nu'\in
\GroupL\cap \GroupA\cap \GroupB^{\pi}$.

\textbf{ b)} Consider sets $\GroupL\cap \GroupA\backslash
\GroupL/\GroupL \cap \GroupB=\{\pi_i^{-1}\,|\,i\in I\}$ and
$\GroupA\cap \GroupL/\GroupA\cap \GroupL\cap
\GroupB^{\pi_i}=\{\nu_{ij}\,|\,j\in J_i\}$ for any $i\in I$. We
should prove that $C=\{\nu_{ij}\pi_i^{-1}\,|\,i\in I,\, j\in
J_i\}$ is a system of representatives $\GroupL/ \GroupL\cap
\GroupB$.

For any $g\in \GroupL$ there is an $i\in I$ such that
$g=a\pi_i^{-1}b$ for some $a\in \GroupL\cap \GroupA$, $b\in
\GroupL\cap \GroupB$. Hence $a=\nu_{ij}a'$ for some $j\in J_i$,
$a'\in \GroupA\cap \GroupL\cap \GroupB^{\pi_i}$. Therefore
$a'=\pi_i^{-1}b'\pi_i$ for some $b'\in \GroupB$. Finally,
$g=\nu_{ij}\pi_i^{-1}b'b$ is contained in $C$ modulo $\GroupL\cap
\GroupB$. Repeating the reasoning from part~a) we obtain that all
elements of $C$ are different modulo $\GroupL\cap \GroupB$.
\end{proof}

\section{Definition and properties of $\DP$}\label{section_F}
In this section we assume that all matrices have entries in a
commutative unitary ring $\Ring$ without divisors of zero.
Moreover, we assume that the characteristic of $\Ring$ is zero.
The case of positive characteristic is discussed in
Remark~\ref{remark1} at the end of this section.

Before defining the function $\DP$ recall the notion of the
pfaffian and the determinant. The pfaffian of a $2r\times 2r$
skew-symmetric matrix $C=(c_{ij})$ is given by %
$$\Pf(C)=\sum\limits_{\{\si\in \Symmgr_{2r}|\, \si(2k-1)<\si(2k)\, \text{ for
all } 1\leq k\leq r\}
 }\sign(\si)\prod\limits_{k=1}^r c_{\si(2k-1),\si(2k)}.$$ %

By the {\it generalized pfaffian} of an arbitrary $2r\times 2r$
matrix $Y=(y_{ij})$ we will mean
$$\P(Y)=\sum\limits_{\si\in \Symmgr_{2r}/\diag(\Symmgr_r\times \Symmgr_r)}
\sign(\si)\prod\limits_{k=1}^r y_{\si(k),\si(k+r)},$$ %
where $\diag(\Symmgr_r\times \Symmgr_r)=\{\nu\times\nu\,|\, \nu\in
\Symmgr_r\}$ is the diagonal subgroup of $\Symmgr_r\times
\Symmgr_r$. By abuse of notation we will also refer to $\P(Y)$ as
the pfaffian. Note that
$$\Pf(Y-Y^{T})=(-1)^{r(r-1)/2}r!\,\P(Y).$$

The determinant of a $t\times t$ matrix  $X=(x_{ij})$ is equal to
$$\det(X)=
\sum\limits_{\si\times\tau\in \Symmgr_{t}\times
\Symmgr_{t}/\diag(\Symmgr_t\times \Symmgr_t)}
\sign(\si\tau)\prod_{k=1}^t x_{\si(k),\tau(k)}.$$

\noindent{\bf Definition 1.} Let  $t,r,s\in\NN$. Given
\begin{itemize}
    \item a $(t+2r)\times (t+2s)$ matrix $X=(x_{ij})$,
    \item a $(t+2r)\times (t+2r)$ matrix $Y=(y_{ij})$,
    \item a $(t+2s)\times (t+2s)$ matrix $Z=(z_{ij})$,
\end{itemize}
we define a function $\DP_{r,s}(X,Y,Z)$,
which is a mixture of the determinant and two pfaffians:%
\begin{eq}\label{eq1}
\sum\limits_{\si\times\tau\in \Symmgr_{t+2r}\times
\Symmgr_{t+2s}/\mathcal{P} }\sign(\si\tau) \prod\limits_{i=1}^t
x_{\si(i),\tau(i)} \prod\limits_{j=1}^r
y_{\si(t+j),\si(t+r+j)}\prod\limits_{k=1}^s
z_{\tau(t+k),\tau(t+s+k)},
\end{eq}
\par \noindent where $\mathcal{P}=
\{\nu_1\times \nu_2\times \nu_2\times \nu_1\times \nu_3\times
\nu_3\,|\,\nu_1\in \Symmgr_t,\nu_2\in \Symmgr_r,\nu_3\in
\Symmgr_s\}$. It is not difficult to see that $\DP_{r,s}(X,Y,Z)$
can also be written in the form
$$
\frac{1}{t!r!s!}\sum\limits_{\si\in \Symmgr_{t+2r}}
\sum\limits_{\tau\in \Symmgr_{t+2s}}\sign(\si\tau)
\prod\limits_{i=1}^t x_{\si(i),\tau(i)} \prod\limits_{j=1}^r
y_{\si(t+j),\si(t+r+j)}\prod\limits_{k=1}^s
z_{\tau(t+k),\tau(t+s+k)}.
$$
In the case $t=r=s=0$ we define $\DP_{r,s}(X,Y,Z)=0$.

\begin{lemma}\label{lemma3}
The function $\DP$ satisfies the following properties:

$1.a)$ $\DP_{r,s}(X^T,Z,Y)=\DP_{r,s}(X,Y,Z)$.

$b)$ $\DP_{r,s}(X,Y^T,Z)=(-1)^r\DP_{r,s}(X,Y,Z)$,

$\quad\,\DP_{r,s}(X,Y,Z^T)=(-1)^s\DP_{r,s}(X,Y,Z)$.

$c)$ if $r=s=0$, then $\DP_{0,0}(X,Y,Z)=\det(X)$; and if $t=0$,
then $\DP_{r,s}(X,Y,Z)=\P(Y)\P(Z)$.

$2$.  Let $g$ be a $(t+2r)\times(t+2r)$ matrix and $h$ be a
$(t+2s)\times(t+2s)$ matrix. Then

$a)$ $\DP_{r,s}(gX,gYg^T,Z)=\det(g)\DP_{r,s}(X,Y,Z)$.

$b)$ $\DP_{r,s}(Xh,Y,h^T Z h)=\det(h)\DP_{r,s}(X,Y,Z)$.

\end{lemma}
\begin{proof} Part 1 follows easily from the definition of $\DP$.

Given mappings $a:[1,t]\to[1,t+2r]$ and $b,c:[1,r]\to[1,t+2r]$,
set

$$M_{a,b,c}= \sum_{\si\in
\Symmgr_{t+2r}} \sum_{\tau\in \Symmgr_{t+2s}} \sign(\si\tau)$$
$$ \prod_{i=1}^t x_{a(i),\tau(i)}
 \prod_{j=1}^r y_{b(j),c(j)}
 \prod_{k=1}^s z_{\tau(t+k),\tau(t+s+k)}
 \prod_{i=1}^t g_{\si(i),a(i)}
 \prod_{j=1}^r g_{\si(t+j),b(j)}  g_{\si(t+r+j),c(j)}.$$

Then
$$\DP_{r,s}(gX,gYg^T,Z)=\frac{1}{t!r!s!}\sum_{a,b,c} M_{a,b,c},$$
where the sum is over all mappings $a:[1,t]\to[1,t+2r]$ and
$b,c:[1,r]\to[1,t+2r]$.

We denote by $(i,j)$ the transposition switching $i$ and $j$ and
make the following observations:

\begin{enumerate}
\item[(i)] If there are $u,v$ with $1\leq u,v\leq t$ such that
$u\neq v$ and $a(u)=a(v)$, then $M_{a,b,c}=0$ because
$\sign(\si)\neq\sign(\si\cdot (u,v))$.

\item[(ii)] If there are $u,v$ with $1\leq u,v\leq r$ such that
$u\neq v$ and $b(u)=b(v)$, then $M_{a,b,c}=0$ because
$\sign(\si)\neq\sign(\si\cdot (t+u,t+v))$.

\item[(iii)] If there are $u,v$ with $1\leq u,v\leq r$ such that
$u\neq v$ and $c(u)=c(v)$, then $M_{a,b,c}=0$ because
$\sign(\si)\neq\sign(\si\cdot (t+r+u,t+r+v))$.

\item[(iv)] If there are $u,v$ with $1\leq u,v\leq r$ such that
$b(u)=c(v)$, then $M_{a,b,c}=0$ because
$\sign(\si)\neq\sign(\si\cdot (t+u,t+r+v))$.

\item[(v)] If there are $u,v$ with $1\leq u\leq t$ and $1\leq
v\leq r$ such that $a(u)=b(v)$, then $M_{a,b,c}=0$ because
$\sign(\si)\neq\sign(\si\cdot (u,t+v))$.

\item[(vi)] If there are $u,v$ with $1\leq u\leq t$ and $1\leq
v\leq r$ such that $a(u)=c(v)$, then $M_{a,b,c}=0$ because
$\sign(\si)\neq\sign(\si\cdot (u,t+r+v))$.
\end{enumerate}

Assume $M_{a,b,c}\neq0$. By (i)--(iii) the mappings $a,b,c$ are
injective and by (iv)--(vi) their images are disjoint. Thus there
is a $\pi\in \Symmgr_{t+2r}$ such that for every $1\leq i\leq t$,
$1\leq j\leq r$ we have $\pi(i)=a(i)$, $\pi(t+j)=b(j)$,
$\pi(t+r+j)=c(j)$. Therefore
$$\DP_{r,s}(gX,gYg^{T},Z)=\frac{1}{t!r!s!} \sum_{\pi\in \Symmgr_{t+2r}}
 \sum_{\si\in \Symmgr_{t+2r}} \sum_{\tau\in \Symmgr_{t+2s}}\sign(\si\tau)$$
$$ \prod_{i=1}^t x_{\pi(i),\tau(i)}
 \prod_{j=1}^r y_{\pi(t+j), \pi(t+r+j)}
 \prod_{k=1}^s z_{\tau(t+k),\tau(t+s+k)}
 \prod_{i=1}^{t+2r} g_{\si\pi^{-1}(i),i}.
$$
Substituting $\si\pi$ for $\si$, we prove 2a). Part 2b) is
analogous.
\end{proof}

Let  $\myun{\ga}\vdash t$, $\myun{\de}\vdash r$, $\myun{\la}\vdash
s$, where $\myun{\ga}=(\my\ga_1,\ldots,\my\ga_u)$,
$\myun{\de}=(\my\de_1,\ldots,\my\de_v)$,
$\myun{\la}=(\my\la_1,\ldots,\my\la_w)$, and let
\begin{itemize}
    \item $X_k=(x_{ij}^k)$ be a $(t+2r)\times(t+2s)$ matrix, where $1\leq k\leq u$,
    \item $Y_k=(y_{ij}^k)$ be a $(t+2r)\times(t+2r)$ matrix, where $1\leq k\leq v$,
  \item $Z_k=(z_{ij}^k)$ be a $(t+2s)\times(t+2s)$ matrix, where $1\leq k\leq w$.
\end{itemize}
Consider the polynomial
$\DP_{r,s}(x_{1}X_1+\ldots+x_{u}X_u,y_1Y_1+\ldots+y_vY_v,z_1Z_1+\ldots+z_wZ_w)$
in variables $x_1,\ldots,x_u$, $y_1,\ldots,y_v$, $z_1,\ldots,z_w$.
Denote by
$$\DP_{\myun{\ga},\myun{\de},\myun{\la}}(X_1,\ldots,X_u,Y_1,\ldots,Y_v,Z_1,\ldots,Z_w)$$ the
coefficient of $x_{1}^{\my\ga_1}\ldots x_{u}^{\my\ga_u}
y_{1}^{\my\de_1}\ldots y_{v}^{\my\de_v} z_{1}^{\my\la_1}\ldots
z_{w}^{\my\la_w}$ in this polynomial. The function
$\DP_{\myun{\ga},\myun{\de},\myun{\la}}$ plays a key role in the
rest of the paper.

\begin{lemma}\label{lemma5}
The polynomial
$\DP_{\myun{\ga},\myun{\de},\myun{\la}}(X_1,\ldots,X_u,Y_1,\ldots,Y_v,Z_1,\ldots,Z_w)$
is equal to
$$\sum\limits_{\si\times \tau\in \Symmgr_{t+2r}\times \Symmgr_{t+2s}/\mathcal{L}}
\sign(\si\tau) \prod\limits_{i=1}^t x^{\My\Ga|i|}_{\si(i),\tau(i)}
\prod\limits_{j=1}^r y^{\My\De|j|}_{\si(t+j),\si(t+r+j)}
\prod\limits_{k=1}^s z^{\My\La|k|}_{\tau(t+k),\tau(t+s+k)},$$
where $\My\Ga,\My\La,\My\De$ are the distributions determined by
$\myun{\ga},\myun{\de},\myun{\la}$, respectively, and
$$\mathcal{L}=\{\nu_1\times\nu_2\times\nu_2\times\nu_1\times\nu_3\times\nu_3\,|\,
\nu_1\in \Symmgr_{\My\Ga},\nu_2\in \Symmgr_{\My\De},\nu_3\in
\Symmgr_{\My\La}\}.$$
\end{lemma}
\begin{proof}
We set
$$\mathcal{A}=\{a:[1,t]\to[1,u],\text{ where }\#a^{-1}(i)=\my\ga_i \text{ for any }1\leq i\leq u\},$$
$$\mathcal{B}=\{b:[1,r]\to[1,v],\text{ where }\#b^{-1}(j)=\my\de_j \text{ for any }1\leq j\leq v\},$$
$$\mathcal{C}=\{c:[1,s]\to[1,w],\text{ where }\#c^{-1}(k)=\my\la_k \text{ for any }1\leq k\leq w\}.$$
It is easy to see that
$$\DP_{\myun{\ga},\myun{\de},\myun{\la}}(X_1,\ldots,X_u,Y_1,\ldots,Y_v,Z_1,\ldots,Z_w)$$
$$=1/t!r!s! \sum_{\si\in \Symmgr_{t+2r}}\sum_{\tau\in \Symmgr_{t+2s}}
\sign(\si\tau) \sum_{a,b,c}
 \prod_{i=1}^t x^{a(i)}_{\si(i),\tau(i)}
 \prod_{j=1}^r y^{b(j)}_{\si(t+j),\si(t+r+j)}
\prod_{k=1}^s z^{c(k)}_{\tau(t+k),\tau(t+s+k)},$$ %
where $a,b,c$ range over $\mathcal{A,B,C}$, respectively. Part~3
of Lemma~\ref{lemma2} gives a bijection between the set
$\Symmgr_t/\Symmgr_{\My\Ga}$ and $\mathcal{A}$ that sends each
$\pi_1\in \Symmgr_t/\Symmgr_{\My\Ga}$ to the mapping
$i\mapsto\My\Ga|\pi_1^{-1}(i)|$. Similarly we obtain bijections
between $\Symmgr_r/\Symmgr_{\My\De}$ and
$\Symmgr_s/\Symmgr_{\My\La}$, respectively, and $\mathcal{B}$ and
$\mathcal{C}$, respectively. Thus
$$\DP_{\myun{\ga},\myun{\de},\myun{\la}}(X_1,\ldots,X_u,Y_1,\ldots,Y_v,Z_1,\ldots,Z_w)=
\frac{1}{t!r!s!}\sum_{\si\in \Symmgr_{t+2r}}\sum_{\tau\in
\Symmgr_{t+2s}}\sign(\si\tau)$$
$$\sum_{\pi_1\in \Symmgr_t/\Symmgr_{\My\Ga}} \sum_{\pi_2\in
\Symmgr_r/\Symmgr_{\My\De}} \sum_{\pi_3\in
\Symmgr_s/\Symmgr_{\My\La}}
  \prod_{i=1}^{t} x^{\My\Ga|\pi_1^{-1}(i)|}_{\si(i),\tau(i)}
  \prod_{j=1}^{r} y^{\My\De|\pi_2^{-1}(j)|}_{\si(t+j),\si(t+r+j)}
\prod_{k=1}^{s}
z^{\My\La|\pi_3^{-1}(k)|}_{\tau(t+k),\tau(t+s+k)}.$$ After
substitution $i\to\pi_1(i), j\to\pi_2(j),k\to\pi_3(k)$ we obtain
the expression
$$\frac{1}{t!r!s!\,\#\Symmgr_{\My\Ga}\, \#\Symmgr_{\My\De}\, \#\Symmgr_{\My\La}}
\sum_{\pi_1\in \Symmgr_t} \sum_{\pi_2\in \Symmgr_r} \sum_{\pi_3\in
\Symmgr_s} \sum_{\si\in \Symmgr_{t+2r}}\sum_{\tau\in
\Symmgr_{t+2s}} \sign(\si\tau)$$
$$      \prod_{i=1}^{t} x^{\My\Ga|i|}_{\si(\pi_1(i)),\tau(\pi_1(i))}
        \prod_{j=1}^{r} y^{\My\De|j|}_{\si(t+\pi_2(j)),\si(t+r+\pi_2(j))}
      \prod_{k=1}^{s}
z^{\My\La|k|}_{\tau(t+\pi_3(k)),\tau(t+s+\pi_3(k))}.$$ Another
substitution $\si\cdot(\pi_1\times\pi_2\times\pi_2)\to\si$,
$\tau\cdot(\pi_1\times\pi_3\times\pi_3)\to\tau$ gives the final
form
$$\frac{1}{\#\Symmgr_{\My\Ga}\,\#\Symmgr_{\My\De}\,\#\Symmgr_{\My\La}}
\sum\limits_{\si\in \Symmgr_{t+2r}}\sum\limits_{\tau\in
\Symmgr_{t+2s}}\sign(\si\tau) \prod\limits_{i=1}^t
x^{\My\Ga|i|}_{\si(i),\tau(i)} \prod\limits_{j=1}^r
y^{\My\De|j|}_{\si(t+j),\si(t+r+j)} \prod\limits_{k=1}^s
z^{\My\La|k|}_{\tau(t+k),\tau(t+s+k)}.
$$
\end{proof}

Denote by $\Ring[x_1,\ldots,x_n]$ and $\QQ[x_1,\ldots,x_n]$,
respectively, the ring of commutative polynomials over $\Ring$ and
$\QQ$, respectively. The following lemma is trivial.

\begin{lemma}\label{lemma_reduction_to_char_0}
Let $p>0$ be the characteristic of $\Ring$. Define the ring
homomorphism $\pi:\ZZ\to\Ring$ by $\pi(1)=1_{\Ring}$, where
$1_{\Ring}$ stands for the identity of $\Ring$. Consider
$a=\sum\al_i a_i\in\Ring[x_1,\ldots,x_n]$, where
$\al_i\in\pi(\Ring)$, and $a_i$ is a monomial in $x_1,\ldots,x_n$.
Take $\be_i\in\{0,1,\ldots,p-1\}\subset\ZZ$ such that
$\pi(\be_i)=\al_i$ and set $b=\sum_i\be_i
a_i\in\QQ[x_1,\ldots,x_n]$.

Then $b=0$ implies $a=0$.
\end{lemma}
\bigskip

We will use the following remark to treat the case of positive
characteristic of $\Ring$.

\begin{remark}\label{remark1}     
When $\Ring$ has a positive characteristic, define $\DP$ by
Formula~\Ref{eq1}. Then Lemmas~\ref{lemma3},~\ref{lemma5} remain
valid.
\end{remark}
\begin{proof}
Analogues of Lemmas~\ref{lemma3} and~\ref{lemma5} assert equality
of certain polynomials over $\pi(\Ring)$ with $\pi$ defined in
Lemma~\ref{lemma_reduction_to_char_0}. Since we have shown that
the corresponding equalities are valid over $\ZZ$, the proof of
Remark~\ref{remark1} follows from
Lemma~\ref{lemma_reduction_to_char_0}.
\end{proof}

\section{Reduction to a zigzag-quiver}\label{section_4}

Throughout this section we use notations from
Section~\ref{section_repr_quiv}.
\bigskip

\noindent{\bf Definition 2}. A quiver $\Quiver$ is called {\it
bipartite}, if every vertex is a source (i.e., there is no arrow
ending at this vertex), or a sink (i.e., there is no arrow starting
at this vertex). A {\it zigzag-quiver} is a bipartite quiver
$\Quiver$ that satisfies the following conditions:
\begin{enumerate}
\item[a)] for each vertex $u\in \Quiver_0$ we have $\phi(u)\neq
u$; moreover, if $u$ is a source, then $\phi(u)$  is a sink and
vice-versa;

\item[b)] there is no arrow $a\in \Quiver_1$ for which
$\al_{a'}=\al_{a''}=\ast$.
\end{enumerate}

Given a zigzag-quiver $\Quiver$, we can depict it schematically as
$$\begin{array}{ccccc}
&\ast\,\bullet     &                   &\bullet\,\,\,\,\,& \\
&\,\,\,\,\vdots  & \longrightarrow   &\vdots\,\,\,\,\,& \\
&\ast\,\bullet     &                   &\bullet\,\,\,\,\,&\\
&                  &  \qquad\nearrow\qquad &                  & \\
&\,\,\,\,\,\bullet &                       &\bullet\,\ast     & \\
&\,\,\,\,\,\vdots  & \longrightarrow       &\vdots\,\,\,\,  & \\
&\,\,\,\,\,\bullet &                        &\bullet\,\ast     &\\
\end{array}
$$ %
where the involution $\phi$ permutes vertices horizontally, and
$\al_u=\ast$ if and only if the vertex $u$  is marked on the
picture by an asterisk. (For more details see the beginning of
Section~\ref{section_results}).


Given a quiver $\Quiver$ and $\mmyun{n},\myun{\al},\phi$,
satisfying requirements~a),~b),~c) from
Section~\ref{section_repr_quiv}, we consider the zigzag-quiver
$\Quiver^{(2)}$ together with
$\mmyun{n}^{(2)},\myun{\al}^{(2)},\phi^{(2)}$, according to the
following construction, which is a variation of ``doubling
construction"{} from~\cite{Schofield92}.
\bigskip

\noindent{\bf Construction 1}. If $\phi(u)=u$, that is a vertex
$u$ does not belong to any pair, then we add a new vertex
$\ov{u}$, which will form a pair with $u$, and change $\mmyun{n}$,
$\phi$, $\myun{\al}$, $SL_{\myun{\al},\phi}(\mmyun{n})$
appropriately. The resulting space of mixed representations of a
new quiver coincides with the space of mixed representations of
the original quiver, because there is no arrow whose head or tail
is $\ov{u}$. So without loss of generality we can assume that each
vertex $u\in \Quiver_0$ belongs to some pair.

Assume $\phi(u)=v\neq u$, where $u,v\in \Quiver_0$ and $\al_u=1$,
$\al_v=\ast$. Add two new vertices $\ov{u}$, $\ov{v}$ to the
quiver and change every arrow $a$ with $a''=u$ to a new arrow
$\ov{a}$ with $\ov{a}''=\ov{u}$ and $\ov{a}'=a'$. Also change
every arrow $a$ with $a''=v$ to a new arrow $\ov{a}$ with
$\ov{a}''=\ov{v}$ and $\ov{a}'=a'$. We will refer to such arrows
as arrows of {\it type} $1$. Define $\al_{\ov{u}}=1$,
$\al_{\ov{v}}=\ast$ and $n_{\ov{u}}=n_{\ov{v}}=n_{u}$. Moreover,
add a new arrow $b$ from $\ov{u}$ to $u$. We will refer to $b$ as
an arrow of {\it type} $2$. Redefine $\phi$ so that
$\phi(u)=\ov{v}$, $\phi(v)=\ov{u}$, and change
$SL_{\myun{\al},\phi}(\mmyun{n})$ appropriately.

After performing this procedure for all pairs, we get the quiver
$\Quiver^{(1)}=(\Quiver_0^{(1)},\Quiver_1^{(1)})$. If an arrow
$a\in \Quiver_1^{(1)}$ satisfies the condition
$\al_{a'}=\al_{a''}=\ast$, change it to a new arrow $c$ such that
$c'=\phi^{(1)}(a'')$, $c''=\phi^{(1)}(a')$. We will refer to $c$
as an arrow of {\it type} $3$. After all of these changes we
obtain a zigzag-quiver
$\Quiver^{(2)}=(\Quiver_0^{(2)},\Quiver_1^{(2)})$ and its space of
mixed representations
$\Reprsp_{\myun{\al}^{(2)},\phi^{(2)}}(\Quiver^{(2)},\mmyun{n}^{(2)})$.
Denote the coordinate functions on
$\Reprsp_{\myun{\al}^{(2)},\phi^{(2)}}(\Quiver^{(2)},\mmyun{n}^{(2)})$
by $y_{ij}^a$ ($a\in \Quiver_1^{(2)}$).

For short we set
$SL^{(1)}=SL_{\myun{\al}^{(1)},\phi^{(1)}}(\mmyun{n}^{(1)})$ and
$SL^{(2)}=SL_{\myun{\al}^{(2)},\phi^{(2)}}(\mmyun{n}^{(2)})$.
\bigskip

\noindent{\bf Example~1}. Let $\Quiver$ be the quiver with vertices
$v,w$, and arrows $a_1,\ldots,a_4$, where $a_1$ goes from $v$ to
$w$, $a_2$ goes in the opposite direction, $a_3$ and $a_4$ are
loops in vertices $v$ and $w$, respectively. Suppose the
involution $\phi$ interchanges vertices $v$ and $w$, $\al_v=1$ and
$\al_w=\ast$. The quiver $\Quiver$ is depicted schematically as
$$a_3\subset \bullet{v} \stackrel{a_1,a_2}{\Longleftrightarrow} \bullet{w} \supset a_4\quad.$$
Then  $\Quiver^{(2)}$ is
$$\begin{array}{ccccccc}
&&\ov{w}\bullet     &  \stackrel{a_2}{\longrightarrow}  &\bullet v&& \\
&&                  &  \qquad\stackrel{a_3,a_4,b}{\nearrow}\qquad &         && \\
&&\ov{v}\bullet &     \stackrel{a_1}{\longrightarrow}   &\bullet w  &,& \\
\end{array}
$$ %
where the arrows of $\Quiver^{(2)}$ are denoted by the same letters as
the corresponding arrows of $\Quiver$, $b$ is a new arrow; the
involution $\phi^{(2)}$ permutes vertices horizontally;  and for a
vertex $i$ of the new quiver $\al^{(2)}_i=\ast$ if and only if
$i=w, \ov{w}$. Note that $a_1,a_2,a_3$ are arrows of type $1$, $b$
is of type $2$, and $a_4$ is of type $3$.

\bigskip
The following theorem shows that for our purpose it is enough to
find out generators for zigzag-quivers.

\begin{theo}\label{theo_reduction}
Consider the homomorphism of $\Field$-algebras %
$$\Phi:
\Field[\Reprsp_{\myun{\al}^{(2)},\phi^{(2)}}(\Quiver^{(2)},\mmyun{n}^{(2)})]
\to
\Field[\Reprsp_{\myun{\al},\phi}(\Quiver,\mmyun{n})],$$ %
defined by
$$\Phi(y_{ij}^a)=\left\{
\begin{array}{rcl}
x_{ij}^a,&&  \text{if the type of }a\text{ is }1 \\
\de_{ij},&&  \text{if the type of }a\text{ is }2 \\
x_{ji}^a,&&  \text{if the type of }a\text{ is }3. \\
\end{array}
\right.
$$
Then its restriction to the algebra of semi-invariants
$$\Field[\Reprsp_{\myun{\al}^{(2)},\phi^{(2)}}
(\Quiver^{(2)},\mmyun{n}^{(2)})]^{SL^{(2)}} \to
\Field[\Reprsp_{\myun{\al},\phi}(\Quiver,\mmyun{n})]^{SL_{\myun{\al},\phi}(\mmyun{n})},$$ %
is a surjective mapping.
\end{theo}
\begin{proof} By Lemma~\ref{lemma1} it is enough to
prove the claim for the quiver $\Quiver^{(1)}$ instead of
$\Quiver^{(2)}$. Moreover, we can  assume that the procedure of
vertex doubling was applied to a single pair of vertices
$\{u,v\}$. Let $b\in \Quiver^{(1)}_1$ be the arrow added to
$\Quiver$ (see the construction of $\Quiver^{(1)}$).

We have $SL^{(1)}=\GroupG_1\times \GroupG_2\times \GroupH$,
$SL_{\myun{\al},\phi}(\mmyun{n})=\GroupG\times \GroupH$ for some
groups $\GroupG,\GroupG_1,\GroupG_2$ that are copies of $SL(n_u)$,
i.e., they are equal to $SL(n_u)$, and some group $\GroupH$. The
group $\GroupG$ acts on $\VectspV_u$ and $\VectspV_v$, the group
$\GroupG_1$ acts on $\VectspV_u^{(1)}$ and
$\VectspV_{\ov{v}}^{(1)}$, and the group $\GroupG_2$ acts on
$\VectspV_{\ov{u}}^{(1)}$ and $\VectspV_{v}^{(1)}$. The group
$\GroupG$ is embedded into $\GroupG_1\times \GroupG_2$ as the
diagonal subgroup
$$\GroupD=\diag(\GroupG_1\times \GroupG_2)=\{(g,g)\,|\,g\in SL(n_u)\}.$$ %
Let $\Quiver^{(1)}\backslash\{b\}$ be a quiver obtained by
removing the arrow $b$ from $\Quiver^{(1)}$ and put
$\Reprsp=\Reprsp_{\myun{\al}^{(1)},\phi^{(1)}}(\Quiver^{(1)}\backslash\{b\},\mmyun{n}^{(1)})$.
The Frobenius reciprocity (\cite{Grosshans97}, Lemma~8.1) implies
that
$$\Psi:\left(\Field[\Reprsp]\otimes \Field\left[SL^{(1)}/{\GroupD\times \GroupH}\right]\right)^{SL^{(1)}}
\to \Field[\Reprsp]^{\GroupD\times \GroupH}
$$  %
is an isomorphism of algebras. Here $SL^{(1)}$ acts on the above
tensor product diagonally and it acts on
${SL^{(1)}}/{\GroupD\times \GroupH}$ by left multiplication. The
mapping $\Psi$ is given by
$$\Psi(f\otimes h)=h(1_{SL^{(1)}})f,$$
where $f\in \Field[\Reprsp]$, $h\in \Field[SL^{(1)}/\GroupD\times
\GroupH]$ and $1_{SL^{(1)}}$ is the identity element of the group
$SL^{(1)}$. On the other hand,
$$\Field[\Reprsp]^{\GroupD\times \GroupH}=\Field[\Reprsp_{\myun{\al},\phi}(\Quiver,\mmyun{n})]^{SL_{\myun{\al},\phi}(\mmyun{n})}.$$

The action of $\GroupG_1\times \GroupG_2$ on $\GroupG_1\times
\GroupG_2/\GroupD$ by left multiplication induces an action of
$SL^{(1)}=\GroupG_1\times \GroupG_2\times \GroupH$ on
$\GroupG_1\times \GroupG_2/\GroupD$. Thus the $SL^{(1)}$-spaces
$SL^{(1)}/\GroupD\times \GroupH$ and $\GroupG_1\times
\GroupG_2/\GroupD$ are isomorphic.  We define a mapping
$\GroupG_1\times \GroupG_2/\GroupD \to SL(n_u)$ by $(g_1,g_2)\to
g_1g_2^{-1}$. Since $(g_1,g_2)\cdot g=g_1gg_2^{-1}$ for $g\in
SL(n_u)$, $g_1\in \GroupG_1$ and $g_2\in \GroupG_2$, this mapping
is an isomorphism of $\GroupG_1\times \GroupG_2$-spaces. Therefore
there is an $SL^{(1)}$-equivariant isomorphism of
$\Field$-algebras
$$\Field[SL^{(1)}/\GroupD\times \GroupH]\simeq \Field[z_{ij}\,|\,1\leq i,j\leq n_u]/I,$$
where $I$ is the ideal generated by the polynomial
$\det((z_{ij})_{1\leq i,j\leq n_u})-1$. Consider the natural
homomorphism
$$\pi: \Field[z_{ij}]\to \Field[z_{ij}]/I.$$ Since
the $SL^{(1)}$-modules $I$ and $\Field[\Reprsp]$ have good
filtrations (see~\cite{Donkin85},~\cite{DZ01}), the mapping $\pi$
induces
a surjective homomorphism  %
$$\ov{\pi}:(\Field[\Reprsp]\otimes \Field[z_{ij}])^{SL^{(1)}}\to
(\Field[\Reprsp]\otimes \Field[SL^{(1)}/\GroupD\times \GroupH])^{SL^{(1)}}.$$ %
Moreover, the mapping $\Field[\Reprsp]\otimes \Field[z_{ij}]\to
\Field[\Reprsp_{\myun{\al}^{(1)},\phi^{(1)}}(\Quiver^{(1)},\mmyun{n}^{(1)})],$
defined by $z_{ij}\to y_{ij}^b$, is an $SL^{(1)}$-equivariant
isomorphism of $\Field$-algebras. It remains to observe that the
restriction of $\Phi$ to semi-invariants is equal to $\Psi\cdot
\ov{\pi}$.
\end{proof}

\section{Semi-invariants of mixed representations of
zigzag-quivers}\label{section_results}
\subsection{Notations}\label{results_notations}

Given a zigzag-quiver $\Quiver=(\Quiver_0,\Quiver_1)$, consider
its mixed representations of a  fixed dimension vector.

Depict $\Quiver$ together with vector spaces, assigned to
vertices, schematically as
$$\begin{array}{ccccl}
\VectspV_1(1)^{\ast}&\bullet  &           &\bullet&\VectspV_1(1) \\
&\vdots                   & \stackrel{j_{+}}{\longrightarrow}        &\vdots& \\
\VectspV_1(l_1)^{\ast}&\bullet&                         &\bullet&\VectspV_1(l_1)\\
                         &&  \stackrel{i}{\nearrow}   && \\
\VectspV_2(1)&\bullet         &                         &\bullet&\VectspV_2(1)^{\ast} \\
&\vdots                   & \stackrel{k_{-}}{\longrightarrow}         &\vdots& \\
\VectspV_2(l_2)&\bullet       &                          &\bullet&\VectspV_2(l_2)^{\ast}.\\
\end{array}
$$ %
Here
\begin{itemize}
\item $\VectspV_1(i)$, $\VectspV_2(j)$, $\VectspV_1(i)^{\ast}$ and
$\VectspV_2(j)^{\ast}$ ($1\leq i\leq l_1$, $1\leq j\leq l_2$) are
vector spaces assigned to vertices. Denote $\dim
\VectspV_1(i)=n_i$ and $\dim \VectspV_2(j)=m_j$. Consider
$\VectspV_1(i)$, $\VectspV_2(j)$ as spaces of column vectors. Fix
the standard bases for $\VectspV_1(i)$, $\VectspV_2(j)$ and the
dual bases for $\VectspV_1(i)^{\ast}$, $\VectspV_2(j)^{\ast}$.

\item the arrows of $\Quiver$ are labeled by symbols
$i,j_{+},k_{-}$, ($1\leq i\leq d_1$, $1\leq j\leq d_2$, $1\leq
k\leq d_3$). Introduce notations $i',i''$, $j'_{+},j''_{+}$,
$k'_{-},k''_{-}$ such that the tail of an arrow $i$ corresponds to
the vector space
$\VectspV_2(i'')$ and its head corresponds to $\VectspV_1(i')$; %
the tail of an arrow $j_{+}$ corresponds to
$\VectspV_1(j''_{+})^{\ast}$ and its head corresponds to $\VectspV_1(j'_{+})$; %
the tail of an arrow $k_{-}$ corresponds to $\VectspV_2(k''_{-})$
and its head corresponds to $\VectspV_2(k'_{-})^{\ast}$.
Schematically
$$\begin{array}{ccccl}
\VectspV_2(i'')&\bullet&\stackrel{i}{\longrightarrow}&\bullet&\VectspV_1(i')\\
\VectspV_1(j''_{+})^{\ast}&\bullet&\stackrel{j_{+}}{\longrightarrow}&\bullet&\VectspV_1(j'_{+})\\
\VectspV_2(k''_{-})&\bullet&\stackrel{k_{-}}{\longrightarrow}&\bullet&\VectspV_2(k'_{-})^{\ast}.\\
\end{array}
$$%

\item the involution $\phi$ permutes vertices of $\Quiver$
horizontally;

\item for a vertex $i$ of zigzag-quiver, $\al_i=\ast$ if and only
if the dual vector space is assigned to $i$.
\end{itemize}

Denote by $\Dimvec$ the dimension vector, whose entries are $n_i$,
$m_j$, where $1\leq i\leq l_1$, $1\leq j\leq l_2$. Then the space
of mixed representations of the quiver $\Quiver$ of dimension
$\Dimvec$
is %
$$\Reprsp=\Reprsp_{\myun{\al},\phi}(\Quiver,\Dimvec) =\oplus_{i=1}^{d_1}\Field^{n_{i'}\times m_{i''}}\bigoplus %
\oplus_{j=1}^{d_2}\Field^{n_{j'_{+}}\times n_{j''_{+}}}\bigoplus
\oplus_{k=1}^{d_3}\Hom_{\Field}\Field^{m_{k'_{-}}\times m_{k''_{-}}} $$ %
$$ \simeq\oplus_{i=1}^{d_1}{\Hom}_{\Field}(\VectspV_2(i''),\VectspV_1(i'))\bigoplus$$ %
$$\oplus_{j=1}^{d_2}\Hom_{\Field}(\VectspV_1(j''_{+})^{\ast},\VectspV_1(j'_{+}))\bigoplus
\oplus_{k=1}^{d_3}\Hom_{\Field}(\VectspV_2(k''_{-}),\VectspV_2(k'_{-})^{\ast}), $$ %
where the isomorphism is given by the choice of bases. Elements of
$\Reprsp$ are mixed representations of $\Quiver$ and we write them
as $H=(H_{1},\ldots,H_{d_1}$, $H_{1}^{+},\ldots,H_{d_2}^{+}$,
$H_{1}^{-},\ldots,H_{d_3}^{-})$. The group
$$\GroupG=SL_{\myun{\al},\phi}(\Dimvec)=\prod_{i=1}^{l_1}SL(n_i)\times
\prod_{j=1}^{l_2}SL(m_j)$$ acts on $\Reprsp$ and on its coordinate
ring $\Field[\Reprsp]$, since $SL(n_i)$ acts on $\VectspV_1(i)$,
$\VectspV_1(i)^{\ast}$, and $SL(m_j)$ acts on $\VectspV_2(j)$,
$\VectspV_2(j)^{\ast}$ (see Section~\ref{section_repr_quiv}).
Denote by

\begin{itemize}
\item $x_{ij}^k$ ($1\leq k\leq d_1$, $1\leq i\leq n_{k'}$, $1\leq
j\leq m_{ k''}$) the coordinate function on $\Reprsp$ that takes a
representation $H$ to the $(i,j)$-th entry of the matrix $H_{k}$;

\item $y_{ij}^k$ ($1\leq k\leq d_2$, $1\leq i\leq n_{k'_{+}}$,
$1\leq j\leq n_{k''_{+}}$) the coordinate function on $\Reprsp$
that takes $H$ to the $(i,j)$-th entry of $H_{k}^{+}$;

\item  $z_{ij}^k$ ($1\leq k\leq d_3$, $1\leq i\leq m_{k'_{-}}$,
$1\leq j\leq m_{k''_{-}}$) the coordinate function on $\Reprsp$
that takes $H$ to the $(i,j)$-th entry of $H_{k}^{-}$.
\end{itemize}%
Denote by $X_k=(x_{ij}^k)$, $Y_k=(y_{ij}^k)$, $Z_k=(z_{ij}^k)$ the
resulting generic matrices.

For $g=(g_{1},\ldots,g_{l_1},h_{1},\ldots,h_{l_2})\in \GroupG$,
$1\leq i\leq d_1$, $1\leq j\leq d_2$, and $1\leq k\leq d_3$ we have %
$$g\cdot X_i= g_{i'}^{-1} X_i h_{i''},\,
g\cdot Y_j= g_{j'_{+}}^{-1} Y_j (g_{j''_{+}}^{-1})^T,\,
g\cdot Z_k= h_{k'_{-}}^{ T} Z_k h_{k''_{-}}.$$ %

\subsection{Formulation of the Theorem}\label{results_general}

Fix $\mmyun{t}=(t_1,\ldots,t_{d_1})$,
$\mmyun{r}=(r_1,\ldots,r_{d_2})$, $\mmyun{s}=(s_1,\ldots,s_{d_3})$
and denote by $\Field[\Reprsp](\mmyun{t},\mmyun{r},\mmyun{s})$ the
space of polynomials that have a total degree $t_i$ in variables
from $X_i$ ($1\leq i\leq d_1$), a total degree $r_j$ in variables
from $Y_j$ ($1\leq j\leq d_2$), and a total degree $s_k$ in
variables from $Z_k$ ($1\leq k \leq d_3$). Further, set
$t=t_1+\ldots +t_{d_1}$, $r=r_1+\ldots+r_{d_2}$,
$s=s_1+\ldots+s_{d_3}$, and let $\MMy{T},\MMy{R},\MMy{S}$ be the
distributions determined by $\mmyun{t},\mmyun{r},\mmyun{s}$,
respectively.
\bigskip

\noindent{\bf  Definition 3.}{\bf (i)} A triplet of multidegrees
$(\mmyun{t},\mmyun{r},\mmyun{s})$ is called {\it admissible} if
there exist $\mmyun{p}=(p_1,\ldots,p_{l_1})\in \NN^{l_1}$,
$\mmyun{q}=(q_1,\ldots,q_{l_2})\in \NN^{l_2}$ such that  %
$$\sum_{1\leq k\leq d_1,\,   k'=i}t_k + \sum_{1\leq k\leq d_2,\, k'_{+}=i}r_k +
  \sum_{1\leq k\leq d_2,\, k''_{+}=i}r_k = n_ip_i  \quad (1\leq i\leq l_1),$$%
$$\sum_{1\leq k\leq d_1,\,  k''=i}t_k + \sum_{1\leq k\leq d_3,\, k'_{-}=i}s_k +
  \sum_{1\leq k\leq d_3,\, k''_{-}=i}s_k = m_iq_i  \quad (1\leq i\leq l_2).$$%
Note that in this case $\sum_{i=1}^{l_1}n_ip_i=t+2r$ and
$\sum_{i=1}^{l_2}m_iq_i=t+2s$. Set $p=p_1+\ldots+p_{l_1}$,
$q=q_1+\ldots+q_{l_2}$, and let $\MMy{P},\MMy{Q}$ be the
distributions determined by $\mmyun{p},\mmyun{q}$, respectively.
We say that $(\mmyun{p},-\mmyun{q})$ is an {\it admissible
weight}.
\smallskip

\noindent{\bf (ii)} Given an admissible triplet
$(\mmyun{t},\mmyun{r},\mmyun{s})$, consider a distribution
$A=(A_1,\ldots,A_p)$ of the set $[1,t+2r]$ and a distribution
$B=(B_1,\ldots,B_q)$ of the set $[1,t+2s]$. If
$\#A_i=n_{\MMy{P}|i|}$ ($1\leq i\leq p$),
$\#B_i=m_{\MMy{Q}|i|}$ ($1\leq i\leq q$),%
$$
\bigcup_{1\leq k\leq d_1,\,  k'=i}\MMy{T}_k %
\bigcup_{1\leq k\leq d_2,\, k'_{+}=i}(t+\MMy{R}_k) %
\bigcup_{1\leq k\leq d_2,\,k''_{+}=i}(t+r+\MMy{R}_k)= %
\bigcup_{1\leq j\leq p,\,\MMy{P}|j|=i}A_j \quad %
(1\leq i\leq l_1)$$
and %
$$
\bigcup_{1\leq k\leq d_1,\, k''=i}\MMy{T}_k %
\bigcup_{1\leq k\leq d_3,\, k'_{-}=i}(t+\MMy{S}_k) %
\bigcup_{1\leq k\leq d_3,\,k''_{-}=i}(t+s+\MMy{S}_k)= %
\bigcup_{1\leq j\leq q,\,\MMy{Q}|j|=i}B_j \quad %
(1\leq i\leq l_2),$$ %
then we say that the quintuple
$(\mmyun{t},\mmyun{r},\mmyun{s},A,B)$ is {\it admissible}.

Given a distribution $A=(A_1,\ldots,A_p)$ of the set $[1,t+2r]$,
define a distribution $A'$ of the set $[1,t]$ and distributions
$A'',A'''$ of the set $[1,r]$ as follows:
$$\begin{array}{cl}
A'=(A_1',\ldots,A_p'),& {\rm where\,\,} A_j'=[1,t] \cap A_j, \\
A''=(A_1'',\ldots,A_p''),& {\rm where\,\,} A_j''=[1,r] \cap (A_j-t),\\
A'''=(A_1''',\ldots,A_p'''),& {\rm where\,\,} A_j''=[1,r] \cap (A_j-t-r).\\
\end{array}$$%
Any empty sets should be omitted. In the same way define
distributions $B',B'',B'''$ for a distribution
$B=(B_1,\ldots,B_q)$ of the set $[1,t+2s]$.
\smallskip

\noindent{\bf (iii)} Let $(\mmyun{t},\mmyun{r},\mmyun{s},A,B)$ be
an admissible quintuple.

There are unique multipartitions
$\myun{\ga}_{max}\vdash\mmyun{t}$, $\myun{\de}_{max}\vdash
\mmyun{r}$, $\myun{\la}_{max}\vdash \mmyun{s}$ such that the
distributions $\My\Ga_{max},\My\De_{max},\My\La_{max}$ determined
by $\myun{\ga}_{max},\myun{\de}_{max},\myun{\la}_{max}$,
respectively, satisfy the equalities:
$$\My\Ga_{max}=A'\cap B'\cap \MMy{T},\,\, %
\My\De_{max}=A''\cap A'''\cap \MMy{R},\,\, %
\My\La_{max}=B''\cap B'''\cap \MMy{S}.$$ Obviously, the octuple
$(\mmyun{t},\mmyun{r},\mmyun{s},A,B,\myun{\ga}_{max},\myun{\de}_{max},\myun{\la}_{max})$
is admissible in the following sense:

Consider a multipartition $\myun{\ga}\vdash \mmyun{t}$ of a
height $\mmyun{u}=(u_1,\ldots,u_{d_1})$, a multipartition
$\myun{\de}\vdash \mmyun{r}$ of a height
$\mmyun{v}=(v_1,\ldots,v_{d_2})$, and a multipartition
$\myun{\la}\vdash \mmyun{s}$ of a height
$\mmyun{w}=(w_1,\ldots,w_{d_3})$ (the definition of {\it height}
was given in Section~\ref{section_Young}). We set
$u=u_1+\ldots+u_{d_1}$, $v=v_1+\ldots+v_{d_2}$,
$w=w_1+\ldots+w_{d_3}$, and we write $\MMy{U},\MMy{V},\MMy{W}$,
respectively, for the distributions determined by
$\mmyun{u},\mmyun{v},\mmyun{w}$, respectively.

If the distributions $\My\Ga,\My\De,\My\La$ determined by
$\myun{\ga},\myun{\de},\myun{\la}$, respectively, satisfy
$$\My\Ga\leq A'\cap B',\,\, \My\De\leq A''\cap A'''\,\,
{\rm and }\,\, \My\La\leq B''\cap B''',$$  then the octuple
$(\mmyun{t},\mmyun{r},\mmyun{s},A,B,\myun{\ga},\myun{\de},\myun{\la})$
is called {\it admissible}. In other words, such an octuple is
admissible if there are mappings
$$a_1:[1,u]\to[1,p],\,\, a_2,a_3:[1,v]\to[1,p],$$
$$b_1:[1,u]\to[1,q],\,\, b_2,b_3:[1,w]\to[1,q]$$
such that
$$\My\Ga_i\subseteq A_{a_1(i)}\cap B_{b_1(i)},$$
$$t+\My\De_j\subseteq A_{a_2(j)},\quad t+r+\My\De_j\subseteq A_{a_3(j)},$$
$$t+\My\La_k\subseteq B_{b_2(k)},\quad t+s+\My\La_k\subseteq B_{b_3(k)}$$
for $1\leq i\leq u$, $1\leq j\leq v$, and $1\leq k\leq w$.

Note that for every $\myun{\ga},\myun{\de},\myun{\la}$ such that
$(\mmyun{t},\mmyun{r},\mmyun{s},A,B,\myun{\ga},\myun{\de},\myun{\la})$
is admissible we have $\myun{\ga}\leq\myun{\ga}_{max}$,
$\myun{\de}\leq\myun{\de}_{max}$,
$\myun{\la}\leq\myun{\la}_{max}$.
\bigskip

\noindent{\bf Construction 2.} Let
$(\mmyun{t},\mmyun{r},\mmyun{s},A,B,\myun{\ga},\myun{\de},\myun{\la})$
be an admissible octuple.
\smallskip

\noindent{\bf (i)}  Given $1\leq k\leq u$, define a $(t+2r)\times
(t+2s)$ matrix $[X]_k$, partitioned into $p\times q$ number of
blocks, where the block in the $(i,j)$-th position is an
$n_{\MMy{P}|i|}\times m_{\MMy{Q}|j|}$ matrix; the block in the
$(a_1(k),b_1(k))$-th position is equal to $X_{\MMy{U}|k|}$, and
the rest of blocks are zero matrices.

Given $1\leq k\leq v$, define a $(t+2r)\times (t+2r)$ matrix
$[Y]_k$, partitioned into  $p\times p$ number of blocks, where the
block in the $(i,j)$-th position is an $n_{\MMy{P}|i|}\times
n_{\MMy{P}|j|}$ matrix; the block in the $(a_2(k),a_3(k))$-th
position is equal to $Y_{\MMy{V}|k|}$, and the rest of blocks are
zero matrices.

Given $1\leq k\leq w$, define a $(t+2s)\times (t+2s)$ matrix
$[Z]_k$, partitioned into $q\times q$ number of blocks, where the
block in the $(i,j)$-th position is a $m_{\MMy{Q}|i|}\times
m_{\MMy{Q}|j|}$ matrix; the block in the $(b_2(k),b_3(k))$-th
position is equal to $Z_{\MMy{W}|k|}$, and the rest of blocks are
zero matrices.
\smallskip

\noindent{\bf (ii)}  For $(\mmyun{t},\mmyun{r},\mmyun{s},A,B)$
consider $\myun{\ga}_{max}$,$\myun{\de}_{max}$,$\myun{\la}_{max}$,
defined in part (iii) of Definition~3. We set
$$\DP_{\myun{\ga},\myun{\de},\myun{\la}}^{A,B}=\DP_{\myun{\ga},\myun{\de},\myun{\la}}(
 [X]_1,\ldots,[X]_u,
 [Y]_1,\ldots,[Y]_v,
 [Z]_1,\ldots,[Z]_w)$$
and
$$\DP_{\mmyun{t},\mmyun{r},\mmyun{s}}^{A,B}=
\DP_{\myun{\ga}_{max},\myun{\de}_{max},\myun{\la}_{max}}^{A,B}.$$
We will see in Lemma~\ref{lemma6} that
$\DP_{\mmyun{t},\mmyun{r},\mmyun{s}}^{A,B}$ is well defined.
\bigskip

\begin{theo}\label{theo1}
If a space
$\Field[\Reprsp](\mmyun{t},\mmyun{r},\mmyun{s})^{\GroupG}$ is
non-zero, then the triplet $(\mmyun{t},\mmyun{r},\mmyun{s})$ is
admissible. In this case
$\Field[\Reprsp](\mmyun{t},\mmyun{r},\mmyun{s})^{\GroupG}$ is
spanned over $\Field$ by the set of all
$\DP_{\mmyun{t},\mmyun{r},\mmyun{s}}^{A,B}$ for admissible
quintuples $(\mmyun{t},\mmyun{r},\mmyun{s},A,B)$.
\end{theo}
\bigskip

Note that in the characteristic zero case Proposition~\ref{prop2}
(see below) provides a formula for calculation of
$\DP_{\mmyun{t},\mmyun{r},\mmyun{s}}^{A,B}\,$, without using
functions $a_1,a_2,a_3,b_1,b_2,b_3$.

Theorem~\ref{theo1} enables us to describe relative invariants
(for definition see Section~\ref{section_repr_quiv}).

\begin{cor}\label{cor}
All relative invariants in $\Field[\Reprsp]$ have admissible
weights. The space of relative invariants of the admissible weight
$(\mmyun{p},-\mmyun{q})$ is spanned over $\Field$ by the set of
all $\DP_{\mmyun{t},\mmyun{r},\mmyun{s}}^{A,B}\,$ for admissible
quintuples $(\mmyun{t},\mmyun{r},\mmyun{s},A,B)$ with weight
$(\mmyun{p},-\mmyun{q})$.
\end{cor}
\begin{proof}See below the proof of Lemma~\ref{lemma6}.
\end{proof}

\subsection{Particular Case of $l_1=l_2=1$}\label{results_special}

To illustrate Theorem~\ref{theo1}, we consider a simple particular
case which contains all the main features of the general case.
Also, during the first reading of the proof of Theorem~\ref{theo1}
it may be useful to work in this particular case; it is quite easy
to rewrite the proof for this case, using the following remarks.

Assume $l_1=l_2=1$. We set $\dim \VectspV_1(1)=n$, $\dim
\VectspV_2(1)=m$. Then Definition~3 and Construction~2 turn into:
\bigskip

\noindent{\bf  Definition 3'.} {(i)} A triplet of multidegrees
$(\mmyun{t},\mmyun{r},\mmyun{s})$ is called {\it admissible} if
there are $p,q\in\NN$ such that
$$t+2r=np,\,\, t+2s=mq.$$

\noindent{(ii)} Let $(\mmyun{t},\mmyun{r},\mmyun{s})$ be an
admissible triplet and $A=(A_1,\ldots,A_p)$, $B=(B_1,\ldots,B_q)$,
respectively, be distributions of the sets $[1,t+2r]$ and
$[1,t+2s]$, respectively.  Then the quintuple
$(\mmyun{t},\mmyun{r},\mmyun{s},A,B)$ is {\it admissible} if
$\#A_i=n$ and $\#B_j=m$ for $1\leq i\leq p$, $1\leq j\leq q$.
\smallskip

\noindent{(iii)} See item (iii) from Definition 3.
\bigskip

\noindent{\bf Construction 2'.} See Construction 2. Note that
sizes of all blocks of matrices $[X]_k$ ($[Y]_k$, $[Z]_k$,
respectively) are $n\times m$ ($n\times n$, $m\times m$,
respectively).
\bigskip

In this case it is not necessary to assume that a quintuple
$(\mmyun{t},\mmyun{r},\mmyun{s},A,B)$ is admissible, hence we
could simplify the result by getting rid of the distributions
$A,B$:

\begin{prop} Let $l_1=l_2=1$. If a triplet $(\mmyun{t},\mmyun{r},\mmyun{s})$ is admissible, then
$\Field[\Reprsp](\mmyun{t},\mmyun{r},\mmyun{s})^{\GroupG}$ is
spanned over $\Field$ by the set of all
$$\DP_{\myun{\ga},\myun{\de},\myun{\la}}(
 [X]_1,\ldots,[X]_u,
 [Y]_1,\ldots,[Y]_v,
 [Z]_1,\ldots,[Z]_w),$$
for arbitrary multipartitions $\myun{\ga}\vdash \mmyun{t}$,
$\myun{\de}\vdash \mmyun{r}$, $\myun{\la}\vdash \mmyun{s}$  of
heights $\mmyun{u}$, $\mmyun{v}$, $\mmyun{w}$, respectively, and
arbitrary mappings
$$a_1:[1,u]\to[1,p],\,\, a_2,a_3:[1,v]\to[1,p],$$
$$b_1:[1,u]\to[1,q],\,\, b_2,b_3:[1,w]\to[1,q],$$
which define $[X]_i$, $[Y]_j$, $[Z]_k$ (see Construction $2$).
\end{prop}
\begin{proof}
The proposition follows immediately from Theorem~\ref{theo1} and
the proof of Lemma~\ref{lemma6} (see below). Note that it is
obvious that matrices $[X]_i$, $[Y]_j$, $[Z]_k$ are well defined.
\end{proof}

\section{Proof of Theorem~\ref{theo1}}\label{section_proof}

\begin{lemma}\label{lemma6}
For an admissible octuple
$(\mmyun{t},\mmyun{r},\mmyun{s},A,B,\myun{\ga},\myun{\de},\myun{\la})$
the element $\DP_{\myun{\ga},\myun{\de},\myun{\la}}^{A,B}$ is well
defined and belongs to
$\Field[\Reprsp](\mmyun{t},\mmyun{r},\mmyun{s})^\GroupG$.
\end{lemma}
\begin{proof}
The first part of the lemma follows by straightforward
calculations.


Given $1\leq k\leq u$ and an $n_{\MMy{P}|a_1(k)|}\times
m_{\MMy{Q}|b_1(k)|}$ matrix $X$, substitute $X$ for the only
non-zero block in $[X]_k$ and denote the resulting
$(t+2r)\times(t+2s)$ matrix by $[X]^{(1)}_k$. In the same manner
for $1\leq k\leq v$ and an $n_{\MMy{P}|a_2(k)|}\times
n_{\MMy{P}|a_3(k)|}$ matrix $Y$ define $(t+2r)\times(t+2r)$ matrix
$[Y]^{(2)}_k$, and for $1\leq k\leq w$ and a
$m_{\MMy{Q}|b_2(k)|}\times m_{\MMy{Q}|b_3(k)|}$ matrix $Z$ define
$(t+2s)\times(t+2s)$ matrix $[Z]^{(3)}_k$.

For  $g=(g_{1},\ldots,g_{l_1},h_{1},\ldots,h_{l_2})\in \GroupG$
define block-diagonal matrices $g_0$, $h_0$ of sizes
$(t+2r)\times(t+2r)$ and $(t+2s)\times(t+2s)$, respectively, in
such a way that the $i$-th block of the matrix $g_{0}$ is equal to
$g_{\MMy{P}|i|}$ ($1\leq i\leq p$) and the $j$-th block of the
matrix $h_{0}$ is equal to $h_{\MMy{Q}|j|}$ ($1\leq j\leq q$). For
variables $x_1,\ldots,x_{u},y_1,\ldots,y_{v},z_1,\ldots,z_{w}$
we have %
$$g\cdot \DP_{r,s}(
\sum_{i=1}^u x_i[X_{\MMy{U}|i|}]_i^{(1)}, \sum_{j=1}^v
y_j[Y_{\MMy{V}|j|}]_j^{(2)},
\sum_{k=1}^w z_k[Z_{\MMy{W}|k|}]_k^{(3)})$$%
$$=\DP_{r,s}(
\sum_{i=1}^u x_i[g_{\MMy{U}|i|'}^{-1} X_{\MMy{U}|i|}\,
h_{\MMy{U}|i|''} ]_i^{(1)},
$$ %
$$
\sum_{j=1}^v y_j[g_{\MMy{V}|j|'_{+}}^{-1}
Y_{\MMy{V}|j|}\,(g_{\MMy{V}|j|''_{+}}^{-1})^T]_j^{(2)},
\sum_{k=1}^w z_k[h_{\MMy{W}|k|'_{-}}^{T}
Z_{\MMy{W}|k|} \,h_{\MMy{W}|k|''_{-}} ]_k^{(3)})$$%
$$=\DP_{r,s}(
g_{0}^{-1}\sum_{i=1}^u x_i[X_{\MMy{U}|i|}]_i^{(1)} h_{0},
g_{0}^{-1}\sum_{j=1}^v y_j[Y_{\MMy{V}|j|}]_j^{(2)}(g_{0}^{-1})^T,
h_{0}^{T}\sum_{k=1}^w z_k[Z_{\MMy{W}|k|}]_k^{(3)} h_{0})$$ which
by Lemma~\ref{lemma3} equals
$$\det(g_{0})\det(h_{0})^{-1} \DP_{r,s}( \sum_{i=1}^u
x_i[X_{\MMy{U}|i|}]_i^{(1)}, \sum_{j=1}^v
y_j[Y_{\MMy{V}|j|}]_j^{(2)}, \sum_{k=1}^w
z_k[Z_{\MMy{W}|k|}]_k^{(3)})$$ and the statement follows.
\end{proof}

We set $$M_1(\mmyun{t},\mmyun{r},\mmyun{s})=
\otimes_{i=1}^{d_1}S^{t_i}(\VectspV_1(i')^{\ast}\otimes
\VectspV_2( i''))\bigotimes $$
$$\otimes_{j=1}^{d_2}S^{r_j}(\VectspV_1(j'_{+})^{\ast}\otimes
\VectspV_1(j''_{+})^{\ast}) \bigotimes
\otimes_{k=1}^{d_3}S^{s_k}(\VectspV_2(k'_{-}) \otimes
\VectspV_2(k''_{-})).$$ By Lemma~\ref{lemma1}, there is an
isomorphism of $\GroupG$-modules
$\Phi:\Field[\Reprsp](\mmyun{t},\mmyun{r},\mmyun{s})\to M_1(\mmyun{t},\mmyun{r},\mmyun{s})$. %
Put
$$M_2(\mmyun{t},\mmyun{r},\mmyun{s},\myun{\ga},\myun{\de},\myun{\la})=
\otimes_{i=1}^{d_1}(\wedge^{\myun{\ga}_i}\VectspV_1(
i')^{\ast}\otimes\wedge^{\myun{\ga}_i}\VectspV_2( i''))\bigotimes
$$
$$\otimes_{j=1}^{d_2}(\wedge^{\myun{\de}_j}\VectspV_1(j'_{+})^{\ast}\otimes
\wedge^{\myun{\de}_j}\VectspV_1(j''_{+})^{\ast})\bigotimes
\otimes_{k=1}^{d_3}(\wedge^{\myun{\la}_k}\VectspV_2(k'_{-})
\otimes\wedge^{\myun{\la}_k}\VectspV_2(k''_{-}))$$ and define a
mapping
$$\zeta_{\myun{\ga},\myun{\de},\myun{\la}}=
\otimes_{i=1}^{d_1}\zeta_{\myun{\ga}_i}\bigotimes
\otimes_{j=1}^{d_2}\zeta_{\myun{\de}_j}\bigotimes
\otimes_{k=1}^{d_3}\zeta_{\myun{\la}_k}:
M_2(\mmyun{t},\mmyun{r},\mmyun{s},\myun{\ga},\myun{\de},\myun{\la})\to
M_1(\mmyun{t},\mmyun{r},\mmyun{s}).$$

\begin{lemma}\label{lemma10}
Given multidegrees $\mmyun{t},\mmyun{r},\mmyun{s}$, we have

$$\Field[\Reprsp](\mmyun{t},\mmyun{r},\mmyun{s})^\GroupG= \sum_{\myun{\ga}\vdash\mmyun{t}}
\sum_{\myun{\de}\vdash\mmyun{r}}
\sum_{\myun{\la}\vdash\mmyun{s}}\Phi^{-1}
\zeta_{\myun{\ga},\myun{\de},\myun{\la}}(M_2(\mmyun{t},\mmyun{r},\mmyun{s},\myun{\ga},\myun{\de},\myun{\la})^\GroupG).$$
\end{lemma}
\begin{proof} Word-by-word repeat the proof
of Proposition~$5.1$ from~\cite{DZ01}. The only difference is that
we consider $M_1(\mmyun{t},\mmyun{r},\mmyun{s})$ instead of
another tensor product from the mentioned Proposition.
\end{proof}

\begin{lemma}\label{lemma9}
(\cite{DZ01}, Proposition~$5.2$) Let
$\myun{\ga}=(\my\ga_1,\ldots,\my\ga_u)\in\NN^u$,
$\my\ga_1+\ldots+\my\ga_u= t$, and let $\VectspV$ be a vector
space with a basis $e_1,\ldots,e_n$. If
$(\wedge^{\myun{\ga}}\VectspV)^{SL(\VectspV)}\neq0$, then $t=np$
for some $p\in\NN$. In this case
$(\wedge^{\myun{\ga}}\VectspV)^{SL(\VectspV)}$ is spanned over
$\Field$ by the elements $\sum_{\tau\in \Symmgr_C/\Symmgr_C\cap
\Symmgr_{\My\Ga}}\sign(\tau)
\otimes_{i=1}^u\wedge_{k\in\My\Ga_i}e_{C\LA\tau(k)\RA}$, where
$C=(C_1,\ldots,C_p)$ is a distribution of $[1,t]$, $\#C_i=n$ for
any $1\leq i\leq p$.
\end{lemma}
\bigskip

Let a quintuple $(\mmyun{t},\mmyun{r},\mmyun{s},A,B)$ be
admissible. Given $\rho_1\in \Symmgr_{t+2r}$ and $\rho_2\in
\Symmgr_{t+2s}$ put
$$F^{A,B}(\rho_1,\rho_2)=\sign(\rho_1\rho_2)\cdot$$ %
\begin{equation}\label{def_of_F}
\prod_{i=1}^t x^{\MMy{T}|i|}_{A\LA \rho_1(i)\RA, B\LA
\rho_2(i)\RA} \prod_{j=1}^r y^{\MMy{R}|j|}_{A\LA \rho_1(t+j)\RA,
A\LA \rho_1(t+r+j)\RA} \prod_{k=1}^s z^{\MMy{S}|k|}_{B\LA
\rho_2(t+k)\RA, B\LA \rho_2(t+s+k)\RA}.
\end{equation}%

The following lemma follows immediately by index change in
products.

\begin{lemma}\label{lemma8}
Let $\rho_1\in \Symmgr_{t+2r}$ and $\rho_2\in \Symmgr_{t+2s}$.
Then

1) $F^{A,B}(\rho_1\cdot (\pi\times id_r\times
id_r),\rho_2)=F^{A,B}(\rho_1,\rho_2\cdot (\pi^{-1}\times
id_s\times id_s))$ for $\pi\in \Symmgr_{\MMy{T}}$.

2) $F^{A,B}(\rho_1\cdot(id_t\times \pi\times
id_r),\rho_2)=F^{A,B}(\rho_1\cdot (id_t\times id_r\times
\pi^{-1}),\rho_2)$ for $\pi\in \Symmgr_{\MMy{R}}$.

3) $F^{A,B}(\rho_1,\rho_2 \cdot (id_t\times \pi\times
id_s))=F^{A,B}(\rho_1,\rho_2 \cdot (id_t\times id_s\times
\pi^{-1}))$ for $\pi\in \Symmgr_{\MMy{S}}$.
\end{lemma}
\bigskip

Given an admissible quintuple
$(\mmyun{t},\mmyun{r},\mmyun{s},A,B)$ and distributions
$\My\Ga,\My\De,\My\La$, respectively, determined by
multipartitions $\myun{\ga}\vdash\mmyun{t}$,
$\myun{\de}\vdash\mmyun{r}$, $\myun{\la}\vdash\mmyun{s}$,
respectively, define
$$H_{\myun{\ga},\myun{\de},\myun{\la}}^{A,B}=
\sum_{\tau_1\in \Symmgr_{A}/\mathcal{L}_1} \sum_{\tau_2\in
\Symmgr_{B}/\mathcal{L}_2} \sum_{\si_1\in
\Symmgr_{\My\Ga}}\sum_{\si_2\in \Symmgr_{\My\De}}\sum_{\si_3\in
\Symmgr_{\My\La}} F^{A,B}(\tau_1\cdot(id_t\times id_r\times
\si_2),\tau_2\cdot(\si_1\times id_s\times
\si_3)),$$%
where $\mathcal{L}_1=\Symmgr_A\cap (\Symmgr_{\My\Ga}\times
\Symmgr_{\My\De}\times \Symmgr_{\My\De})$,
$\mathcal{L}_2=\Symmgr_B\cap (\Symmgr_{\My\Ga}\times
\Symmgr_{\My\La}\times \Symmgr_{\My\La})$.

\begin{prop}\label{prop1}
1. If a triplet $(\mmyun{t},\mmyun{r},\mmyun{s})$ is not
admissible, then
$\Field[\Reprsp](\mmyun{t},\mmyun{r},\mmyun{s})^\GroupG=0$.

2. If a triplet $(\mmyun{t},\mmyun{r},\mmyun{s})$ is admissible,
then $\Field[\Reprsp](\mmyun{t},\mmyun{r},\mmyun{s})^\GroupG$ is
spanned over $\Field$ by the set
$\{H_{\myun{\ga},\myun{\de},\myun{\la}}^{A,B}\,|$ a quintuple
$(\mmyun{t},\mmyun{r},\mmyun{s},A,B)$ is admissible,
$\myun{\ga}\vdash\mmyun{t}$, $\myun{\de}\vdash\mmyun{r}$,
$\myun{\la}\vdash\mmyun{s}\}$.
\end{prop}
\begin{proof}
Rewrite
$M_2(\mmyun{t},\mmyun{r},\mmyun{s},\myun{\ga},\myun{\de},\myun{\la})^\GroupG$
as
$$
\begin{array}{l}
\ot_{i=1}^{l_1} \Big( \ot\limits_{
\begin{subarray}{c}0< k\leq d_1\\k'=i\end{subarray}
} \wedge^{\myun{\ga}_k}\VectspV_1(k')^{\ast}
\bigotimes %
\ot\limits_{
\begin{subarray}{c}0< k\leq d_2\\k'_{+}=i\end{subarray}
} \wedge^{\myun{\de}_k}\VectspV_1(
k'_{+})^{\ast} \bigotimes %
\ot\limits_{
\begin{subarray}{c}0< k\leq d_2\\k''_{+}=i\end{subarray}
}
\wedge^{\myun{\de}_k}\VectspV_1(k''_{+})^{\ast} %
\Big)^{SL(n_i)}
\bigotimes\\ %
\ot_{j=1}^{l_2}\Big( \ot\limits_{
\begin{subarray}{c}0< k\leq d_1\\k''=j\end{subarray}
} \wedge^{\myun{\ga}_k}\VectspV_2( k'')\bigotimes %
\ot\limits_{
\begin{subarray}{c}0< k\leq d_3\\k'_{-}=j\end{subarray}
}\wedge^{\myun{\la}_k}\VectspV_2(k'_{-})\bigotimes %
\ot\limits_{
\begin{subarray}{c}0< k\leq d_3\\k''_{-}=j\end{subarray}
} \wedge^{\myun{\la}_k}\VectspV_2(k''_{-})\Big)^{SL(m_j)}
\end{array}$$
and use Lemmas~\ref{lemma10} and~\ref{lemma9}.
\end{proof}

\begin{remark}\label{remark2}
From now on we work over the field $\QQ$ instead of a field
$\Field$. The case of positive characteristic is considered at the
end of the section (see Remark~\ref{remark_end_of_section}).
We will use formulas for
$\DP_{\myun{\ga},\myun{\de},\myun{\la}}^{A,B}$ and
$H_{\myun{\ga},\myun{\de},\myun{\la}}^{A,B}$ with denominators,
i.e.,
$$H_{\myun{\ga},\myun{\de},\myun{\la}}^{A,B}= \frac{1}{c}
\sum_{\tau_1\in \Symmgr_A}\sum_{\tau_2\in \Symmgr_B}
\sum_{\si_1\in \Symmgr_{\My\Ga}}\sum_{\si_2\in
\Symmgr_{\My\De}}\sum_{\si_3\in \Symmgr_{\My\La}}
F^{A,B}(\tau_1\cdot(id_t,id_r,\si_2),\tau_2\cdot(\si_1,id_s,\si_3)),$$%
where $c= \#(\Symmgr_A\cap (\Symmgr_{\My\Ga}\times
\Symmgr_{\My\De}\times \Symmgr_{\My\De}))\, \#(\Symmgr_B\cap
(\Symmgr_{\My\Ga}\times \Symmgr_{\My\La}\times
\Symmgr_{\My\La}))$, and similarly for
$\DP_{\myun{\ga},\myun{\de},\myun{\la}}^{A,B}$.
\end{remark}

\begin{prop}\label{prop2}
Assume an octuple
$(\mmyun{t},\mmyun{r},\mmyun{s},A,B,\myun{\ga},\myun{\de},\myun{\la})$
is admissible. Consider permutations $\pi_1\in \Symmgr_{t+2r}$,
$\pi_2\in \Symmgr_{t+2s}$ such that $A=\MMy{N}^{\pi_1}$,
$B=\MMy{M}^{\pi_2}$ and $A\LA i\RA=\MMy{N}\LA\pi_1(i)\RA$, $B\LA
j\RA=\MMy{M}\LA\pi_2(j)\RA$ for $1\leq i\leq t+2r$, $1\leq j\leq
t+2s$, where $\MMy{N},\MMy{M}$ are distributions determined by the
vectors
$\mmyun{n}=(n_1,\ldots,n_1,\ldots,n_{l_1},\ldots,n_{l_1})\in\NN^{p_1+\ldots+p_{l_1}}$,
$\mmyun{m}=(m_1,\ldots,m_1,\ldots,m_{l_2},\ldots,m_{l_2})\in\NN^{q_1+\ldots+q_{l_2}}$,
respectively.
Then %
$$\DP_{\myun{\ga},\myun{\de},\myun{\la}}^{A,B}=\sign(\pi_1\pi_2)
\frac{1}{\#\Symmgr_{\My\Ga}\,\#\Symmgr_{\My\De}\,\#\Symmgr_{\My\La}}
\sum_{\tau_1\in \Symmgr_A}\sum_{\tau_2\in \Symmgr_B}
F^{A,B}(\tau_1,\tau_2),
$$%
where $F^{A,B}$ is defined by formula~\Ref{def_of_F}.
\end{prop}
\begin{proof}
Note that the admissibility of the octuple implies
$(\Symmgr_{\My\Ga}\times \Symmgr_{\My\De}\times
\Symmgr_{\My\De})\subseteq \Symmgr_A$, $(\Symmgr_{\My\Ga}\times
\Symmgr_{\My\La}\times \Symmgr_{\My\La})\subseteq \Symmgr_B$.

Define $[X]_k=(f_{ij}^k)$, $[Y]_k=(g_{ij}^k)$ and
$[Z]_k=(h_{ij}^k)$.
By Lemma~\ref{lemma5} and the definition%
$$
\DP_{\myun{\ga},\myun{\de},\myun{\la}}^{A,B}=
\frac{1}{\#\Symmgr_{\My\Ga}\,\#\Symmgr_{\My\De}\,\#\Symmgr_{\My\La}}
\sum_{\tau_1\in \Symmgr_{t+2r}}\sum_{\tau_2\in
\Symmgr_{t+2s}}\sign(\tau_1\tau_2) L,\mbox{ where }$$
$$L=
\prod_{i=1}^t f^{\My\Ga|i|}_{\tau_1(i),\tau_2(i)} \prod_{j=1}^r
g^{\My\De|j|}_{\tau_1(t+j),\tau_1(t+r+j)} \prod_{k=1}^s
h^{\My\La|k|}_{\tau_2(t+k),\tau_2(t+s+k)}.
$$%
For any $1\leq i\leq t$, $1\leq j\leq r$, and $1\leq k\leq s$ we have %
$$
f^{\My\Ga|i|}_{\tau_1(i),\tau_2(i)}= \left\{
\begin{array}{cl}
x^{\MMy{U}|\My\Ga|i||}_{\MMy{N}\LA\tau_1(i)\RA,\MMy{M}\LA\tau_2(i)\RA},&
\mbox{ if } \tau_1(i)\in\MMy{N}_{a_1(\My\Ga|i|)} \, \mbox{and} \, \\
&\quad\,\,\tau_2(i)\in\MMy{M}_{b_1(\My\Ga|i|)}\\
0,&\mbox{ otherwise },\\
\end{array}
\right.
$$%
$$
g^{\My\De|j|}_{\tau_1(t+j),\tau_1(t+r+j)}= \left\{
\begin{array}{cl}
y^{\MMy{V}|\My\De|j||}_{\MMy{N}\LA\tau_1(t+j)\RA,\MMy{N}\LA\tau_1(t+r+j)\RA},&
\mbox{ if }
\tau_1(t+j)  \in\MMy{N}_{a_2(\My\De|j|)}\, \mbox{and} \, \\
&\quad\,\, \tau_1(t+r+j)\in\MMy{N}_{a_3(\My\De|j|)}\\
0,&\mbox{ otherwise },\\
\end{array}
\right.
$$%
$$
h^{\My\La|k|}_{\tau_2(t+k),\tau_2(t+s+k)}= \left\{
\begin{array}{cl}
z^{\MMy{W}|\My\La|k||}_{\MMy{M}\LA\tau_2(t+k)\RA,\MMy{M}\LA\tau_2(t+s+k)\RA},&
\mbox{ if }
\tau_2(t+k)  \in\MMy{M}_{b_2(\My\La|k|)}\, \mbox{and} \,\\
&\quad\,\,\tau_2(t+s+k)\in\MMy{M}_{b_3(\My\La|k|)}\\
0,&\mbox{ otherwise }.\\
\end{array}
\right.
$$%
Since  the octuple is admissible, it implies
$i\in\MMy{N}^{\pi_1}_{a_1(\My\Ga|i|)}\cap
\MMy{M}^{\pi_2}_{b_1(\My\Ga|i|)},$
$$t+j  \in \MMy{N}^{\pi_1}_{a_2(\My\De|j|)},\,\,
t+r+j\in \MMy{N}^{\pi_1}_{a_3(\My\De|j|)},$$
$$t+k  \in\MMy{M}^{\pi_2}_{b_2(\My\La|k|)},\,\,
t+s+k\in \MMy{M}^{\pi_2}_{b_3(\My\La|k|)}.$$ %
Therefore
$$\pi_1(i)\in \MMy{N}_{a_1(\My\Ga|i|)},\,\,
\pi_1(t+j)\in \MMy{N}_{a_2(\My\De|j|)},\,\,
\pi_1(t+r+j)\in \MMy{N}_{a_3(\My\De|j|)},$$ %
$$\pi_2(i)\in \MMy{M}_{b_1(\My\Ga|i|)},\,\,
\pi_2(t+k)\in \MMy{M}_{b_2(\My\La|k|)},\,\,
\pi_2(t+s+k)\in \MMy{M}_{b_3(\My\La|k|)}.$$ %
If $\tau_1\in \Symmgr_{\MMy{N}}\pi_1$ and $\tau_2\in
\Symmgr_{\MMy{M}}\pi_2$, then $L$ equals %
$$
\prod\limits_{i=1}^t
x^{\MMy{T}|i|}_{\MMy{N}\LA\tau_1(i)\RA,\MMy{M}\LA\tau_2(i)\RA}
\prod\limits_{j=1}^r
y^{\MMy{R}|j|}_{\MMy{N}\LA\tau_1(t+j)\RA,\MMy{N}\LA\tau_1(t+r+j)\RA}
\prod\limits_{k=1}^s
z^{\MMy{S}|k|}_{\MMy{M}\LA\tau_2(t+k)\RA,\MMy{M}\LA\tau_2(t+s+k)\RA};
$$%
otherwise $L$ is equal to zero. Use a substitution
$\tau_1\to\pi_1\si_1$, $\tau_2\to\pi_2\si_2$ in the given
expression for $\DP_{\myun{\ga},\myun{\de},\myun{\la}}^{A,B}$ to
conclude the proof.
\end{proof}

\begin{cor}\label{cor1}  
Suppose an octuple
$(\mmyun{t},\mmyun{r},\mmyun{s},A,B,\myun{\ga},\myun{\de},\myun{\la})$
is admissible and $\myun{\ga}_{max}$, $\myun{\de}_{max}$,
$\myun{\la}_{max}$ corresponds to
$(\mmyun{t},\mmyun{r},\mmyun{s},A,B)$. Then
$$\DP_{\myun{\ga},\myun{\de},\myun{\la}}^{A,B}=\pm c\,
\DP_{\myun{\ga}_{max},\myun{\de}_{max},\myun{\la}_{max}}^{A,B},$$ %
where
$c=\#\Symmgr_{\My\Ga_{max}}\#\Symmgr_{\My\De_{max}}\#\Symmgr_{\My\La_{max}}\Big/
\#\Symmgr_{\My\Ga}\,\#\Symmgr_{\My\De}\,\#\Symmgr_{\My\La}\in\NN$.
\end{cor}

\begin{lemma}\label{lemma11}  
Assume a quintuple $(\mmyun{t},\mmyun{r},\mmyun{s},A,B)$ is admissible,
$\myun{\ga}\vdash\mmyun{t}$, $\myun{\de}\vdash\mmyun{r}$, and
$\myun{\la}\vdash\mmyun{s}$. Then there are $\myun{\ga}_0\vdash\mmyun{t}$,
$\myun{\de}_0\vdash\mmyun{r}$, $\myun{\la}_0\vdash\mmyun{s}$, $\rho_1\in
\Symmgr_{\My\Ga}$, $\rho_2\in \Symmgr_{\My\De}$, and $\rho_3\in \Symmgr_{\My\La}$ such
that the octuple $(\mmyun{t},\mmyun{r},\mmyun{s},A^{\rho_1\times\rho_2\times\rho_2},
B^{\rho_1\times\rho_3\times\rho_3}, \myun{\ga}_0,\myun{\de}_0,\myun{\la}_0)$ is
admissible, $\Symmgr_{\My\Ga_0}=\Symmgr_{A'}^{\rho_1}\cap \Symmgr_{\My\Ga}\cap
\Symmgr_{B'}^{\rho_1}$, $\Symmgr_{\My\De_0}=\Symmgr_{A''}^{\rho_2}\cap
\Symmgr_{\My\De}\cap \Symmgr_{A'''}^{\rho_2}$, and
$\Symmgr_{\My\La_0}=\Symmgr_{B''}^{\rho_3}\cap \Symmgr_{\My\La}\cap
\Symmgr_{B'''}^{\rho_3}$, where $\My\Ga,\My\De,\My\La,\My\Ga_0,\My\De_0,\My\La_0$ are the
distributions determined by
$\myun{\ga},\myun{\de},\myun{\la},\myun{\ga}_0,\myun{\de}_0,\myun{\la}_0$, respectively.
\end{lemma}
\begin{proof}
Intersecting  the distributions $A',\My\Ga,B'$ of the set $[1,t]$,
we obtain the distribution $\My\Ga_0^{\rho_1^{-1}}$ for some
$\myun{\ga}_0\vdash\mmyun{t}$, $\rho_1\in \Symmgr_{\My\Ga}$.
Intersecting the distributions $A'',\My\De,A'''$ of the set
$[1,r]$, we obtain the distribution $\My\De_0^{\rho_2^{-1}}$ for
some $\myun{\de}_0\vdash\mmyun{r}$, $\rho_2\in \Symmgr_{\My\De}$.
Intersecting the distributions $B'',\My\La,B'''$ of the set
$[1,s]$, we obtain the distribution $\My\La_0^{\rho_3^{-1}}$ for
some $\myun{\la}_0\vdash\mmyun{s}$, $\rho_3\in \Symmgr_{\My\La}$.
Such elements
$\myun{\ga}_0,\myun{\de}_0,\myun{\la}_0,\rho_1,\rho_2,\rho_3$
satisfy the required properties.
\end{proof}

\begin{prop}\label{prop3}
Assume that a quintuple ($\mmyun{t},\mmyun{r},\mmyun{s},A,B$) is
admissible, $\myun{\ga}\vdash\mmyun{t}$,
$\myun{\de}\vdash\mmyun{r}$, $\myun{\la}\vdash\mmyun{s}$, and
$\My\Ga,\My\De,\My\La$ are distributions determined by
$\myun{\ga},\myun{\de},\myun{\la}$, respectively. Then
the semi-invariant $H_{\myun{\ga},\myun{\de},\myun{\la}}^{A,B}$ is equal to %
$$
\sum_{\pi_1^{-1}\in \Symmgr_{\My\Ga}\cap \Symmgr_{A'} \backslash
\Symmgr_{\My\Ga}/\Symmgr_{\My\Ga}\cap \Symmgr_{B'}}\;
\sum_{\pi_2^{-1}\in \Symmgr_{\My\De}\cap \Symmgr_{A''}\backslash
\Symmgr_{\My\De}/\Symmgr_{\My\De}\cap \Symmgr_{A'''}}\;
\sum_{\pi_3^{-1}\in \Symmgr_{\My\La}\cap \Symmgr_{B''}\backslash
\Symmgr_{\My\La}/\Symmgr_{\My\La}\cap \Symmgr_{B'''}}
\!\!\!\!\!\!\pm \DP_{\myun{\ga}_0,\myun{\de}_0,\myun{\la}_0}^{
A_0,B_0},
$$%
where $\myun{\ga}_0\vdash\mmyun{t}$, $\myun{\de}_0\vdash\mmyun{r}$,
$\myun{\la}_0\vdash\mmyun{s}$, $A_0=A^{\rho_1\times\rho_2\times\pi_2\rho_2}$,
$B_0=B^{\pi_1\rho_1\times\rho_3\times\pi_3\rho_3}$, and $\rho_1\in \Symmgr_t$, $\rho_2\in
\Symmgr_r$, $\rho_3\in \Symmgr_s$ depend on
$\pi_1,\pi_2,\pi_3,\mmyun{t},\mmyun{r},\mmyun{s},A,B,\myun{\ga},\myun{\de},\myun{\la}$.
Moreover, the octuple $(\mmyun{t},\mmyun{r},\mmyun{s}, A_0, B_0,
\myun{\ga}_0,\myun{\de}_0,\myun{\la}_0)$ is admissible.
\end{prop}
\begin{proof}
Since $\Symmgr_{\My\Ga}\subseteq \Symmgr_{\MMy{T}}$,
$\Symmgr_{\My\De}\subseteq \Symmgr_{\MMy{R}}$ and
$\Symmgr_{\My\La}\subseteq \Symmgr_{\MMy{S}}$ we can use
Lemma~\ref{lemma8} throughout the proof. Denote %
$$c_1=  %
\#((\Symmgr_{\My\Ga}\times \Symmgr_{\My\De}\times \Symmgr_{\My\De})\cap %
\Symmgr_{A})\,\, \#((\Symmgr_{\My\Ga}\times \Symmgr_{\My\La}\times  %
\Symmgr_{\My\La})\cap \Symmgr_{B})$$  %
$$=\#(\Symmgr_{\My\Ga}\cap \Symmgr_{A'})\,\cdot\, %
\#(\Symmgr_{\My\De}\cap \Symmgr_{A''})\,\cdot\, %
\#(\Symmgr_{\My\De}\cap \Symmgr_{A'''})\,\cdot\, %
\#(\Symmgr_{\My\Ga}\cap \Symmgr_{B'})\,\cdot\, %
\#(\Symmgr_{\My\La}\cap \Symmgr_{B''})\,\cdot\, %
\#(\Symmgr_{\My\La}\cap \Symmgr_{B'''}).$$ %
For $1\leq i\leq 3$ substitute $\si_i^{-1}$ for $\si_i$ in the
definition
of $H_{\myun{\ga},\myun{\de},\myun{\la}}^{A,B}$ to get %
$$H_{\myun{\ga},\myun{\de},\myun{\la}}^{A,B}=
\frac{1}{c_1} \sum_{\tau_1\in \Symmgr_A, \tau_2\in \Symmgr_B}
\sum_{\si_1\in \Symmgr_{\My\Ga}} \sum_{\si_2\in \Symmgr_{\My\De}}
\sum_{\si_3\in \Symmgr_{\My\La}} F^{A,B}(\tau_1\cdot (id_t\times
id_r\times \si_2^{-1}),\tau_2\cdot
(\si_1^{-1}\times id_s\times \si_3^{-1})).$$%
Define $\mathcal{M}_1=\Symmgr_{\My\Ga}/\Symmgr_{\My\Ga}\cap
\Symmgr_{B'},\,\,
\mathcal{M}_2=\Symmgr_{\My\De}/\Symmgr_{\My\De}\cap
\Symmgr_{A'''},\,\,
\mathcal{M}_3=\Symmgr_{\My\La}/\Symmgr_{\My\La}\cap
\Symmgr_{B'''},$
$$\mathcal{L}_1=\Symmgr_{\My\Ga}\cap \Symmgr_{B'}, \,\,
\mathcal{L}_2=\Symmgr_{\My\De}\cap \Symmgr_{A'''}, \,\,
\mathcal{L}_3=\Symmgr_{\My\La}\cap \Symmgr_{B'''}.$$ Then
$H_{\myun{\ga},\myun{\de},\myun{\la}}^{A,B}$
$$=\frac{1}{c_1}
\sum_{\substack{\si_i\in \mathcal{M}_i\\1\leq i\leq 3}}\;
\sum_{\substack{\nu_i\in \mathcal{L}_i\\1\leq i\leq 3}}\;
\sum_{\substack{\tau_1\in \Symmgr_A\\\tau_2\in \Symmgr_B}}\;
\!\!\!\! F^{A,B}( \tau_1\cdot(id_t\times id_r\times
\nu_2^{-1}\si_2^{-1}), \tau_2\cdot (\nu_1^{-1}\si_1^{-1}\times
id_s\times \nu_3^{-1}\si_3^{-1}))$$%
which is equal to
$$
\frac{1}{c_2} \sum_{\si_i\in \mathcal{M}_i,1\leq i\leq 3}
\sum_{\tau_1\in \Symmgr_A,\tau_2\in \Symmgr_B}F^{A,B}( \tau_1
\cdot(\si_1\times \si_2\times id_r), \tau_2 \cdot (id_t\times
\si_3\times id_s)),
$$%
where $c_2= \#(\Symmgr_{\My\Ga}\cap \Symmgr_{A'})\,
\#(\Symmgr_{\My\De}\cap \Symmgr_{A''})\, \#(\Symmgr_{\My\La}\cap
\Symmgr_{B''})$, after we substituted $\tau_1\cdot (id_t\times
id_r\times \nu_2^{-1})\to\tau_1$ and $\tau_2\cdot
(\nu_1^{-1}\times id_s\times \nu_3^{-1})\to\tau_2$. Define
$$\mathcal{\Field}_1=\Symmgr_{\My\Ga}\cap \Symmgr_{A'}\backslash \Symmgr_{\My\Ga}/
\Symmgr_{\My\Ga}\cap \Symmgr_{B'}, $$
$$\mathcal{\Field}_2=\Symmgr_{\My\De}\cap \Symmgr_{A''}\backslash
\Symmgr_{\My\De}/ \Symmgr_{\My\De}\cap \Symmgr_{A'''}, \,\,
\mathcal{\Field}_3=\Symmgr_{\My\La}\cap \Symmgr_{B''}\backslash
\Symmgr_{\My\La}/
\Symmgr_{\My\La}\cap \Symmgr_{B'''}, $$ %
$$\mathcal{T}_1=\Symmgr_{A' }\cap \Symmgr_{\My\Ga}/\Symmgr_{A'} \cap \Symmgr_{\My\Ga}\cap
\Symmgr_{B'  }^{\pi_1}, $$
$$\mathcal{T}_2=\Symmgr_{A''}\cap
\Symmgr_{\My\De}/\Symmgr_{A''}\cap \Symmgr_{\My\De}\cap
\Symmgr_{A'''}^{\pi_2}, \,\, \mathcal{T}_3=\Symmgr_{B''}\cap
\Symmgr_{\My\La}/\Symmgr_{B''}\cap \Symmgr_{\My\La}\cap
\Symmgr_{B'''}^{\pi_3}. $$ %
Lemma~\ref{lemma4} implies $\mathcal{M}_i=\{\nu_i\pi_i^{-1}\,|\,
\pi_i^{-1}\in \mathcal{\Field}_i,\,\nu_i\in \mathcal{T}_i\}$ for
$1\leq i\leq 3$.
Thus $H_{\myun{\ga},\myun{\de},\myun{\la}}^{A,B}$ equals %
$$
\frac{1}{c_2} %
\sum_{\substack{\pi_i^{-1}\in\mathcal{\Field}_i\\1\leq i\leq 3}}\; %
\sum_{\substack{\nu_i     \in \mathcal{T}_i\\1\leq i\leq 3}}\;
\sum_{\substack{\tau_1\in \Symmgr_A\\\tau_2\in \Symmgr_B}}
F^{A,B}( \tau_1\cdot (\nu_1\pi_1^{-1}\times \nu_2\pi_2^{-1}\times
id_r), \tau_2\cdot (id_t\times \nu_3\pi_3^{-1}\times id_s)).
$$%
Substitute $\tau_1$ for  $\tau_1\cdot (\nu_1\times \nu_2\times
id_r)$ and $\tau_2$ for $\tau_2\cdot (id_t\times \nu_3\times
id_s)$ and define
$$c_3= %
\#(\Symmgr_{A'} \cap \Symmgr_{\My\Ga}\cap \Symmgr_{B'  }^{\pi_1})\,\cdot\, %
\#(\Symmgr_{A''}\cap \Symmgr_{\My\De}\cap \Symmgr_{A'''}^{\pi_2})\,\cdot\, %
\#(\Symmgr_{B''}\cap \Symmgr_{\My\La}\cap
\Symmgr_{B'''}^{\pi_3}).$$
Then %
$$H_{\myun{\ga},\myun{\de},\myun{\la}}^{A,B}=
\sum_{\pi_i^{-1}\in \mathcal{\Field}_i,1\leq i\leq 3}
\sum_{\tau_1\in \Symmgr_A,\tau_2\in \Symmgr_B} \frac{1}{c_3}
F^{A,B}(\tau_1\cdot(id_t\times id_r\times
\pi_2),\tau_2\cdot(\pi_1\times id_s\times \pi_3))
$$%
$$=\sum_{\pi_i^{-1}\in \mathcal{\Field}_i,1\leq i\leq 3}
\frac{1}{c_3} \sum_{\tau_1\in \Symmgr_A^{id_t\times id_r\times
\pi_2}}\sum_{\tau_2\in \Symmgr_B^{\pi_1\times id_s\times \pi_3}}
F^{A,B}((id_t\times id_r\times \pi_2)\cdot\tau_1,(\pi_1\times
id_s\times \pi_3)\cdot\tau_2).
$$%
Applying Lemma~\ref{lemma11} to an admissible quintuple
$(\mmyun{t},\mmyun{r},\mmyun{s},A^{id_t\times id_r\times
\pi_2},B^{\pi_1\times id_r\times \pi_3})$ we obtain
$\myun{\ga}_0\vdash\mmyun{t}$, $\myun{\de}_0\vdash\mmyun{r}$,
$\myun{\la}_0\vdash\mmyun{s}$, $\rho_1\in \Symmgr_{\My\Ga}$,
$\rho_2\in \Symmgr_{\My\De}$ and $\rho_3\in \Symmgr_{\My\La}$ such
that the octuple $(\mmyun{t},\mmyun{r},\mmyun{s}, A^{\rho_1\times
\rho_2\times \pi_2\rho_2},B^{\pi_1\rho_1\times \rho_3\times
\pi_3\rho_3}, \myun{\ga}_0,\myun{\de}_0,\myun{\la}_0)$ is
admissible, $\Symmgr_{\My\Ga_0}=\Symmgr_{A' }^{\rho_1}\cap
\Symmgr_{\My\Ga}\cap \Symmgr_{B'}^{\pi_1\rho_1}$,
$\Symmgr_{\My\De_0}=\Symmgr_{A''}^{\rho_2}\cap
\Symmgr_{\My\De}\cap \Symmgr_{A'''}^{\pi_2\rho_2}$ and
$\Symmgr_{\My\La_0}=\Symmgr_{B''}^{\rho_3}\cap
\Symmgr_{\My\La}\cap \Symmgr_{B'''}^{\pi_3\rho_3}$. If we set
$\rho_{I}=\rho_1\times \rho_2\times \rho_2$,
$\rho_{II}=\rho_1\times \rho_3\times \rho_3$, $\pi_{I}=id_t\times
id_r\times \pi_2$, and $\pi_{II}=\pi_1\times id_s\times \pi_3$,
then
$$
H_{\myun{\ga},\myun{\de},\myun{\la}}^{A,B}= \sum_{\pi_i^{-1}\in
\mathcal{\Field}_i,i\in{1,3}} \frac{1}{c_3} \sum_{\tau_1\in
\Symmgr_A^{\pi_{I}\rho_{I}},\tau_2\in
\Symmgr_B^{\pi_{II}\rho_{II}}}
F^{A,B}(\pi_{I}\rho_{I}\tau_1,\pi_{II}\rho_{II}\tau_2).
$$%
By Lemma~\ref{lemma2}, there are permutations $\eta_1\in
\Symmgr_A^{\pi_{I}\rho_{I}}$ and  $\eta_2\in
\Symmgr_B^{\pi_{II}\rho_{II}}$ that satisfy the conditions
$A\LA\pi_{I }\rho_{I }(i)\RA=A^{\pi_{I }\rho_{I
}}\LA\eta_1(i)\RA$,
$B\LA\pi_{II}\rho_{II}(j)\RA=B^{\pi_{II}\rho_{II}}\LA\eta_2(j)\RA$
for $1\leq i\leq t+2r$, $1\leq j\leq t+2s$. Hence
$$F^{A,B}(\pi_{I}\rho_{I}\tau_1,\pi_{II}\rho_{II}\tau_2)=
\sign(\pi_1\pi_2\pi_3\tau_1\tau_2)\prod_{i=1}^t x^{\MMy{T}|i|}_{
A^{\pi_{I }\rho_{I }}\LA \eta_1\tau_1(i)\RA,
B^{\pi_{II}\rho_{II}}\LA \eta_2\tau_2(i)\RA}$$
$$\prod_{j=1}^r
y^{\MMy{R}|j|}_{ A^{\pi_{I}\rho_{I}}\LA \eta_1\tau_1(t+j)\RA,
A^{\pi_{I}\rho_{I}}\LA \eta_1\tau_1(t+r+j)\RA}\prod_{k=1}^s
z^{\MMy{S}|k|}_{ B^{\pi_{II}\rho_{II}}\LA \eta_2\tau_2(t+k)\RA,
B^{\pi_{II}\rho_{II}}\LA \eta_2\tau_2(t+s+k)\RA}.$$ %
Finally, substitute $\tau_1$ for $\eta_1\tau_1$, and  $\tau_2$ for
$\eta_2\tau_2$. Using Proposition~\ref{prop2} and the equality
$c_3=\#\Symmgr_{\My\Ga_0}\,\#\Symmgr_{\My\De_0}\,\#\Symmgr_{\My\La_0}$
we derive the required assertion.
\end{proof}

Corollary~\ref{cor1} and Proposition~\ref{prop3} show that for any
admissible quintuple ($\mmyun{t},\mmyun{r},\mmyun{s},A_0,B_0$) the
element $H_{\my\ga,\my\de,\my\la}^{A,B}$ belongs to the span over
$\ZZ$ of the
set %
$$\{\DP_{r,s}^{A,B}\,|\text{ a quintuple }
(\mmyun{t},\mmyun{r},\mmyun{s},A,B)\text{ is admissible}\}.$$ %

\begin{remark}\label{remark_end_of_section}
Using Lemma~\ref{lemma_reduction_to_char_0} together with
reasoning from the proof of Remark~\ref{remark1} it is easy to see
that the last statement remains valid over a field $\Field$ of
positive characteristic.
\end{remark}
Proposition~\ref{prop1} concludes the proof of
Theorem~\ref{theo1}.


\section{An example}\label{section_appl}

We will apply Theorem~\ref{theo1} to the following example.

Let $\Quiver$ be  a zigzag-quiver with $l_1=0$, $l_2=1$, $m=2$
(see Section~\ref{section_results}). For short, we write $d$ for
$d_3$. Depict $\Quiver$ schematically as
$$\begin{array}{ccccc}
\Field^2&\bullet &  \stackrel{1_{-},\ldots,d_{-}}{\longrightarrow} &\bullet&(\Field^2)^{\ast} \\
\end{array} $$%
The space of mixed representations is $\Reprsp=(\Field^{2\times
2})^d$, where, as usual, $\Field^{2\times 2}$ is the space of all
$2\times 2$ matrices over $\Field$. The group $\GroupG=SL(2)$ acts
on the coordinate ring $\Field[\Reprsp]=\Field[z_{ij}^k\,|\,1\leq
i,j\leq 2,\,1\leq k\leq d]$ by the formula: $g\cdot Z_k=g^T Z_k
g$, $g\in \GroupG$. Note that $\Reprsp$ corresponds to $d$-tuple
of bilinear forms on $\Field^2$ (see Section~\ref{section_intro}
for details). Invariants of pairs of bilinear forms were
investigated in~\cite{Adamovich_Golovina77}.

Consider the action of $GL(2)$ on $\Field[\Reprsp]$ by $g\cdot
Z_k=g^{-1} Z_k g$, $g\in GL(2)$. Let
$$J=
\left(
\begin{array}{cc}
0 & 1 \\
-1 & 0 \\
\end{array}
\right)
$$
Define an automorphism of algebras  $\Phi: \Field[\Reprsp]\to
\Field[\Reprsp]$ such that $\Phi(z_{ij}^k)$ is equal to the
$(i,j)$-th entry of $Z_k J$.

\begin{prop}\label{prop_appl}
The restriction of $\Phi$ to $GL(2)$-invariants is an isomorphism
of algebras
$$\Field[\Reprsp]^{GL(2)}\cong \Field[\Reprsp]^\GroupG.$$
\end{prop}
\begin{proof}
The algebra $\Field[\Reprsp]^{GL(2)}$ is known to be generated by
$\det(Z_k)$,  $\tr(Z_{k_1}\ldots Z_{k_r})$, where $1\leq
k,k_1,\ldots,k_r\leq d$ (see~\cite{Donkin92a}). Denote by
$\GroupA(d)$ the $\Field$-algebra generated by $\det(Z_k)$,
$\tr(Z_{k_1}J\ldots Z_{k_r}J)$. We have to prove
$\Field[\Reprsp]^{\GroupG}=\GroupA(d)$.

Note that for $g\in\GroupG$ we have $gJg^T=J$ and %
$$g\cdot
\tr(Z_{k_1}J\cdot Z_{k_2}J\cdot \ldots\cdot Z_{k_r}J)$$
$$=\tr(g^T Z_{k_1}gJ\cdot g^T Z_{k_2}gJ\cdot \ldots \cdot
g^T Z_{k_r}gJ)=\tr(Z_{k_1}J\cdot Z_{k_2}J\cdot \ldots \cdot
Z_{k_r}J).$$ %
This proves the inclusion $\GroupA(d) \subset
\Field[\Reprsp]^{\GroupG}$.

Theorem~\ref{theo1} implies that $\Field[\Reprsp]^\GroupG$ is
generated by $\DP^{0,B}_{0,0,\mmyun{s}}$, where
$B=(B_1,\ldots,B_q)$, $\mmyun{s}=(s_1,\ldots,s_d)$, and a
quintuple $(0,0,\mmyun{s},0,B)$ is admissible. So
$q=s=s_1+\ldots+s_d$.  Our goal is to show  that
$\DP^{0,B}_{0,0,\mmyun{s}}$ is a polynomial with integer
coefficients in the above mentioned generators of $\GroupA(d)$.
Obviously, without loss of generality we can assume $\Field=\QQ$.


Let $\myun{\la}_{max}\vdash \mmyun{s}$ be the partition given in
part (iii) of Definition~3, and let $\MMy{S},\My\La_{max}$ be the
distributions determined by $\mmyun{s},\myun{\la}_{max}$,
respectively. Since $\My\La_{max}=B''\cap B'''\cap \MMy{S}$, the
following formula holds:
$$\#\Symmgr_{\My\La_{max}}=\prod_{k=1}^d  \prod_{1\leq i,j\leq s}
\#\{l\in[1,s]\,:\,B|l|=i,\,B|l+s|=j,\,\MMy{S}|l|=k\}!.$$ %
We denote the right hand side of the above equality by
$c(B,\mmyun{s})$.

For an arbitrary $\mmyun{s}\in\NN^d$ and a distribution
$B=(B_1,\ldots,B_s)$ of $[1,2s]$ with $\#B_1=\ldots=\#B_s=2$
define
$$\P^{B}_{\mmyun{s}}=\frac{1}{c(B,\mmyun{s})}\sum_{\tau\in \Symmgr_B} \prod_{l=1}^s
\sign(\tau) z_{B\LA \tau(l)\RA,B\LA \tau(s+l)\RA}^{\MMy{S}|l|}.$$
By Proposition~\ref{prop2}, we have
$\DP^{0,B}_{0,0,\mmyun{s}}=\pm\P^{B}_{\mmyun{s}}$  for the
admissible quintuple $(0,0,\mmyun{s},0,B)$.

A distribution $B$ is called {\it decomposable}, if there is a
proper subset $\Quiver$ of $[1,s]$ such that for $1\leq l\leq s$
we have $B|l|,\, B|s+l|\in \Quiver$, or $B|l|,\, B|s+l|\notin
\Quiver$. If $B$ is decomposable, then there are $\mmyun{s}_r$,
$B_r$ and a mapping $\psi_r: [1,d]\to [1,d]$ (where $r$ is $1$ or
$2$) that satisfy
$\P_{\mmyun{s}}^{B}=\ov{\psi}_1(\P_{\mmyun{s}_1}^{B_1})
\ov{\psi}_2(\P_{\mmyun{s}_2}^{B_2})$.  Here the homomorphism of
algebras  $\ov{\psi}_r:\Field[\Reprsp]\to \Field[\Reprsp]$ is
given by substitution $\ov{\psi}_r(z_{ij}^k)=z_{ij}^{\psi_r(k)}$.
Therefore without loss of generality we can assume that $B$ is
indecomposable.

The indecomposability of $B$ implies that $c(B,\mmyun{s})$ is 1 or
2. Moreover, $c(B,\mmyun{s})= 2$ if and only if $s=2$,
$B=(\{1,2\},\{3,4\})$, and in that case
$\P^{B}_{\mmyun{s}}=\det(Z_1)$. Therefore we can assume
$c(B,\mmyun{s})= 1$.

If $s=1$, then $B=(\{1,2\})$ and
$\P^{B}_{\mmyun{s}}=z_{12}^1-z_{21}^1=-\tr(Z_1J)$. Hence we can
assume $s>1$. Let $D$ be a distribution of $[1,2s]$ defined by:
$D_1=\{1,2s\}$, $D_l=\{l,s+l-1\}$, where $2\leq l\leq s-1$,
$D_s=\{s,2s-1\}$. Then
$$\P_{\mmyun{s}}^D=\sum_{\tau_1,\ldots,\tau_s\in \Symmgr_2}
\sign(\tau_1\ldots\tau_s) z_{\tau_1(1),\tau_2(2)}^{\MMy{S}|1|}
z_{\tau_2(1),\tau_3(2)}^{\MMy{S}|2|} \ldots
z_{\tau_s(1),\tau_1(2)}^{\MMy{S}|s|}.$$

For ${\be}=(\be_1,\ldots,\be_d)$ with $\be_k$ equals 1 or 2,
define homomorphism of $\Field$-algebras
$\psi_{{\be}}:\Field[\Reprsp]\to \Field[\Reprsp]$ by
$$\psi_{{\be}}(z_{ij}^k)=
\left\{
\begin{array}{cl}
z_{ij}^k, &\text{ if } \be_k=1 \\
z_{ji}^k, &\text{ if } \be_k=2. \\
\end{array}
\right.$$ It is not difficult to see that there is a $\be$ and a
mapping $\psi:[1,d]\to[1,d]$ such that
$\P^{B}_{\mmyun{s}}=\pm\ov{\psi}(\psi_{{\be}}(\P^{D}_{\mmyun{s}}))$.

For $1\leq k\leq d$ define $U_k=Z_kJ$ and denote by $u_{ij}^k$ the
entries of $U_k$. Hence $u_{ij}^k=(-1)^j z_{i,\xi(j)}^k$, where
$\xi=(12)\in \Symmgr_2$ is a transposition. We obtain
$$\tr(U_{\MMy{S}|1|}\ldots U_{\MMy{S}|s|})=\sum_{1\leq i_1,\ldots,i_s \leq 2}
u_{i_1i_2}^{\MMy{S}|1|} u_{i_2 i_3}^{\MMy{S}|2|}\ldots u_{i_s
i_1}^{\MMy{S}|s|}.$$ For $1\leq l\leq s$ define a permutation
$\tau_l$ on $\{1,2\}$ by $\tau_l(1)=i_l$. Thus
$\tau_l(2)=\xi(i_l)$ and $\sign(\tau_l)=-(-1)^{i_l}$. This implies
$$\tr(U_{\MMy{S}|1|}\ldots U_{\MMy{S}|s|})=(-1)^s \P_{\mmyun{s}}^{D}.$$

Summarizing, for an indecomposable distribution $B$ we have shown
that $\P_{\mmyun{s}}^{B}=\pm\tr(V_{k_1}J\ldots V_{k_s}J)$, where
$V_k$ is $Z_k$ or $Z_k^{T}$. The formulas
$$\tr(H_1^TJM_2)=\tr(H_1JM_2)-\tr(H_1J)\tr(H_2),\,\,\,\,\tr(H_1^TJ)=-\tr(H_1J),$$
valid for $2\times 2$ matrices $H_1,H_2$, conclude the proof.
\end{proof}

\begin{remark}
Proposition~\ref{prop_appl} can not be proven by elementary
methods, because for all $g,h\in SL(2)$ the inequality
$g^{-1}Z_1Jg\neq h^T Z_1 h$ holds. Moreover, if $g,h\in SL(2)$ and
$\GroupB\subset \Field^{2\times 2}$ consists of elements $H$ that
satisfy $g^{-1} H J g= h^T H h$, then the dimension of the closure
of $\GroupB$ in Zariski topology is strictly less than $4$.
\end{remark}

\bigskip
\noindent{\bf Acknowledgements.} The authors were supported by
``Complex integration projects of SB RAS 1.9"{} and by RFFI (grant
04.01.00489). The authors are grateful to the referees whose
comments improved the paper considerably.

\end{document}